\documentclass[11pt]{article}

\usepackage{geometry}
\usepackage{morefloats}
\usepackage{amsfonts}
\usepackage{color}
\usepackage{multirow}
\usepackage{latexsym}
\usepackage{amsmath,amssymb,amsthm}
\usepackage{dsfont}
\usepackage{extarrows}
\usepackage[authoryear, sort]{natbib}
\usepackage[hidelinks=true]{hyperref}

\numberwithin{equation}{section}

\newtheoremstyle{thm}
{9pt}
{9pt}
{\itshape}
{}
{\bfseries}
{.}
{ }
{}
\theoremstyle{thm}

\newtheorem{theorem}{Theorem}[section]
\newtheorem{lemma}[theorem]{Lemma}
\newtheorem{corollary}[theorem]{Corollary}

\newtheoremstyle{def}
{9pt}
{9pt}
{}
{}
{\bfseries}
{.}
{ }
{}
\theoremstyle{def}

\newtheorem{definition}[theorem]{Definition}
\newtheorem{remark}[theorem]{Remark}
\newtheorem{example}[theorem]{Example}

\newcommand{\R}{\mathbb{R}} 
\newcommand{\Z}{\mathbb{Z}} 
\newcommand{\E}{\mathbb{E}} 

\renewcommand{\footnoterule}{%
	\kern -3.5pt
	\hrule width \textwidth height 1pt
	\kern 3.5pt
}

\makeatletter
\def\blfootnote{\xdef\@thefnmark{}\@footnotetext}
\makeatother

\allowdisplaybreaks

\begin{document}

\title{\bf Fixed point characterizations of continuous univariate probability distributions and their applications}


\author{\Large Steffen Betsch and Bruno Ebner \\ \normalsize Karlsruhe Institute of Technology (KIT), Institute of Stochastics \\ \normalsize Karlsruhe, Germany \\ \\}

\date{\today}
\maketitle


\begin{abstract}
	By extrapolating the explicit formula of the zero-bias distribution occurring in the context of Stein's method, we construct characterization identities for a large class of absolutely continuous univariate distributions. Instead of trying to derive characterizing distributional transformations that inherit certain structures for the use in further theoretic endeavours, we focus on explicit representations given through a formula for the density- or distribution function. The results we establish with this ambition feature immediate applications in the area of goodness-of-fit testing. We draw up a blueprint for the construction of tests of fit that include procedures for many distributions for which little (if any) practicable tests are known. To illustrate this last point, we construct a test for the Burr Type XII distribution for which, to our knowledge, not a single test is known aside from the classical universal procedures.
\end{abstract}

\newpage

\section{Introduction}
\label{SEC Intro}
Over the last decades, Stein's method for distributional approximation has become a viable tool for proving limit theorems and establishing convergence rates. At it's heart lies the well-known Stein characterization which states that a real-valued random variable $Z$ has a standard normal distribution if, and only if,
\begin{equation} \label{Stein char normal distr}
	\E \big[ f^{\prime}(Z) - Z f(Z) \big] = 0
\end{equation}
holds for all functions $f$ of a sufficiently large class of test functions. To exploit this characterization for testing the hypothesis
\begin{equation}\label{Hyp0}
	\mathbf{H_0}: \mathbb{P}^X \in \big\{ \mathcal{N}(\mu, \sigma^2) \, | \, (\mu, \sigma^2) \in \R \times (0,\infty) \big\}
\end{equation}
of normality, where $\mathbb{P}^X$ is the distribution of a real-valued random variable $X$, against general alternatives, \cite{BE:2019:1} used that (\ref{Stein char normal distr}) can be untied from the class of test functions with the help of the so-called zero-bias transformation introduced by \cite{GR:1997}. To be specific, a real-valued random variable $X^*$ is said to have the $X$-zero-bias distribution if
\begin{equation*}
	\E \big[ f^{\prime}(X^*) \big] = \E \big[ X f(X) \big]
\end{equation*}
holds for any of the respective test functions $f$. If $\E X = 0$ and $\mathrm{Var}(X)=1$, the $X$-zero-bias distribution exists and is unique, and it has distribution function
\begin{equation} \label{explicit formula zero-bias trafo}
	T^X (t) = \E \big[ X \, (X - t) \, \mathds{1}\{ X \leq t \} \big], \quad t \in \R .
\end{equation}
By (\ref{Stein char normal distr}), the standard Gaussian distribution is the unique fixed point of the transformation $\mathbb{P}^X \mapsto \mathbb{P}^{X^*}$. Thus, the distribution of $X$ is standard normal if, and only if,
\begin{equation} \label{fixed point statement normal distr}
	T^X = F_X ,
\end{equation}
where $F_X$ denotes the distribution function of $X$. In the spirit of characterization-based goodness-of-fit tests, an idea introduced by \cite{L:1953}, this fixed point property directly admits a new class of testing procedures as follows. Letting $\widehat{T}^X_n$ be an empirical version of $T^X$ and $\widehat{F}_n$ the empirical distribution function, both based on the standardized sample, \cite{BE:2019:1} proposed a test for (\ref{Hyp0}) based on the statistic
\begin{equation*}
	G_n = n \int_{\R} \left| \widehat{T}^X_n (t) - \widehat{F}_n (t) \right|^2 w(t) \, \mathrm{d}t ,
\end{equation*}
where $w$ is an appropriate weight function, which, in view of (\ref{fixed point statement normal distr}), rejects the normality hypothesis for large values of $G_n$. As these tests have several desirable properties such as consistency against general alternatives, and since they show a very promising performance in simulations, we devote this work to the question to what extent the fixed point property and the class of goodness-of-fit procedures may be generalized to other distributions. 

Naturally, interest in applying Stein's method to other distributions has already grown and delivered some corresponding results. Characterizations like (\ref{Stein char normal distr}) have been established en mass [for an overview on characterizing Stein operators and further references, we recommend the work by \cite{LRS:2017}]. Charles Stein himself presented some ideas fundamental to the so-called density approach [see \cite{S:1986}, Chapter VI, and \cite{SDHR:2004}, Section 5] which we shall use as the basis of our considerations. Related results for the special case of exponential families were already given by \cite{H:1978} and \cite{R:1979}. Another approach pioneered by \cite{B:1990} [see also \cite{G:1991}] includes working with the generator of the semi-group of operators corresponding to a Markov process whose stationary distribution is the one in consideration. A third advance is based on fixed point properties of probability transformations like the zero-bias transformation. Very general distributional transformations were introduced by \cite{GR:2005} and refined by \cite{D:2017}. In the latter contribution, the transformations, and with them the explicit formulae, rely heavily on sign changes of the so-called biasing functions. These sign changes, in turn, depend on the parameters of the distribution in consideration which renders the explicit representations impractical for the use in goodness-of-fit testing.

The starting point of the present paper is the density approach identity. Here, a result more general than (\ref{Stein char normal distr}) is provided by showing that, for suitable density functions $p$, a given real-valued random variable $X$ has density $p$ if, and only if,
\begin{equation} \label{Stein chara density approach}
	\E \left[ f^{\prime}(X) + \frac{p^{\prime}(X)}{p(X)} f(X) \right] = 0
\end{equation}
holds for a sufficiently large class of test functions. We provide fixed point characterizations like (\ref{fixed point statement normal distr}) by using the analogy between (\ref{Stein chara density approach}) and (\ref{Stein char normal distr}) to extrapolate the explicit formula (\ref{explicit formula zero-bias trafo}) of the zero-bias transformation to other distributions. Using this approach, these transformations will no longer be probability transformations, but we maintain the characterizing identity which suffices for the use in goodness-of-fit testing. Our confidence in the approach is manifested by the fact that it has already been implemented by \cite{BE:2019:2} for the special case of the Gamma distribution.

With our results we contribute to the growing amount of applications of Stein's- (or the Stein-Chen) method and his characterization in the realm of statistics. Much has been done in the area of stochastic modeling, which often includes statistical methods. For instance, \cite{F:2014} and \cite{RR:2010} [see also \cite{B:1982} and \cite{B:1989}] tackle counting problems in the context of random graphs with Stein's method. The technique led to further insights in time series- and mean field analysis, cf. \cite{K:2000} and \cite{Y:2017}. \cite{BD:2017} and \cite{BDF:2016} developed Stein's method for diffusion approximation which is used as a tool for performance analysis in the theory of queues. As for statistical research that is more relatable to our pursuits, quite a bit is known when it comes to normal approximation for maximum likelihood estimators, investigated, for instance, by \cite{A:2018}, \cite{AG:2018}, \cite{AR:2017}, and \cite{P:2017}, to name but a few contributions. Moreover, \cite{GPR:2017} consider chi-square approximation to study Pearson's statistic which is used for goodness-of-fit testing in classification problems. Also note that \cite{AR:2018} apply the results of \cite{GPR:2017} to obtain bounds to the chi-square distribution for twice the log-likelihood ratio, the statistic used for the classical likelihood ratio test. Finally, the contributions by \cite{CSG:2016} and \cite{LLJ:2016} aim at the goal we also pursue in Section \ref{goodness-of-fit}, namely to apply Steinian characterizations to construct goodness-of-fit tests for probability distributions.

The paper at hand is organized as follows. We first introduce an appropriate setting for our considerations by stating the conditions for a density function to fit into our framework and prove identity (\ref{Stein chara density approach}) in this specific setting. We then give our characterization results, distinguishing between distributions supported by the whole real line, those with semi-bounded support and distributions with bounded support. Throughout, we give examples of density functions of different nature to show that our conditions are not restrictive, as well as to provide connections to characterizations that are already known and included in our statements. Next, we consider applications in goodness-of-fit testing and show that the proposed tests include the classical Kolmogorov-Smirnov- and Cram\'{e}r-von Mises procedures, as well as three modern tests considered in the literature. To illustrate the methods in the last part, we construct the first ever goodness-of-fit test specifically for the two parameter Burr Type XII distribution, and show in a simulation study that the test is sound and powerful compared to classical procedures.

\section{Notation and regularity conditions}
\label{SEC Notation}

Throughout, let $(\Omega, \mathcal{A}, \mathbb{P})$ be a probability space and $p$ a non-negative density function supported by an interval $\mathrm{spt}(p) = [L, R]$, where $- \infty \leq L < R \leq \infty$, and with $\int_L^R p(x) \, \mathrm{d}x = 1$. Denoting by $P$ the distribution function associated with $p$, we state the following regularity conditions:
\begin{itemize}
	\item[\bfseries{(C1)}]{The function $p$ is continuous and positive on $(L, R)$, and there are $L < y_1 < \dotso < y_m < R$ such that $p$ is continuously differentiable on $(L, y_1)$, $(y_\ell, y_{\ell + 1})$, $\ell \in \{ 1, \dots, m-1 \}$, and $(y_m , R)$.}
\end{itemize}
Whenever (C1) holds, we write $S(p) = (L, R) \setminus \{ y_1, \dots, y_m \}$.
\begin{itemize}	
	\item[\bfseries{(C2)}]{For the map $S(p) \ni x \mapsto \kappa_p(x) = \left| \frac{p^{\prime}(x) \min\{ P(x), \, 1 - P(x) \}}{p^2(x)} \right|$ we have $\sup_{x \, \in \, S(p)} \kappa_p(x) < \infty$,}
	\item[\bfseries{(C3)}]{$\int_{S(p)} \big( 1 + |x|\big) \big|p^{\prime}(x)\big| \mathrm{d}x < \infty$,}
	\item[\bfseries{(C4)}]{$\lim_{x \, \searrow \, L} \frac{P(x)}{p(x)} = 0$, and}
	\item[\bfseries{(C5)}]{$\lim_{x \, \nearrow \, R} \frac{1 - P(x)}{p(x)} = 0$.}
\end{itemize}
The integral $\int_{S(p)}$ is understood as the sum of the integrals over the interval-components in (C1). For a probability density function $p$ that satisfies (C1), and a function $f : (L, R) \to \R$ which is differentiable on $S(p)$ except in one point, we denote the point of non-differentiability in $S(p)$ by $t_f$, and set
\begin{align*}
	S(p, f)
	= S(p) \setminus \big\{ t_f \big\}
	= (L, R) \setminus \big\{ y_1, \dots, y_m, t_f \big\} .
\end{align*}
We index the elements of $(L, R) \setminus S(p, f)$ with $y_1^f < \dotso < y_{m + 1}^{f}$.
\begin{definition}[\textbf{Test functions}] \label{DEF test functions}
	For a probability density function $p$ with $\mathrm{spt}(p) = [L, R]$ that satisfies (C1), we denote by $\mathcal{F}_p$ the set of all functions $f : (L, R) \to \R$ that are continuous on $(L, R)$ and differentiable on $S(p)$ except in (precisely) one point, that satisfy
	\begin{align*}
	\sideset{}{_{x \, \searrow \, L}}\lim f(x) \, p(x)
	= \sideset{}{_{x \, \nearrow \, R}}\lim f(x) \, p(x)
	= 0,
	\end{align*}
	and for which $x \mapsto \frac{p^{\prime}(x)}{p(x)} f(x)$ and $x \mapsto f^{\prime}(x)$ are bounded on $S(p, f)$.
\end{definition}

We write $\mathcal{L}^1$ for the Borel-Lebesgue measure on the real line, and $X \sim p \mathcal{L}^1$ whenever a random variable $X$ has Lebesgue density $p$, and we define
\begin{align*}
	\mathrm{disc}(X)
	= \big\{ t \in (L, R) \, \big| \, \mathbb{P}(X = t) > 0 \big\},
\end{align*}
the set of all atoms of a random variable $X$, containing at most countably many points.

\section{The density approach identity}
\label{SEC Density approach}

In this section, we restate the density approach identity. Since we use a very particular class of test functions which bears some technicalities, we give an outline of the proof in Appendix \ref{APP proof of density approach}, roughly following \cite{LS:2013}. We refer to Section II of \cite{LS:2013:2} for a discrete version of the density approach identity, and mention \cite{LS:2016} for related statements in the context of parametric distributions.
\begin{lemma} \label{LEMMA density approach}
	If $p$ is a probability density function with $\mathrm{spt}(p) = [L, R]$ that satisfies (C1) and (C2), and if $X : \Omega \to (L, R)$ is an arbitrary random variable with $\mathbb{P}\big( X \in S(p) \big) = 1$, then $X \sim p \mathcal{L}^1$ if, and only if,
	\begin{align*}
	\E \left[ f^{\prime}(X) + \frac{p^{\prime}(X)}{p(X)} \, f(X) \right] = 0
	\end{align*}
	for each $f \in \mathcal{F}_p$ with $t_f \notin \mathrm{disc}(X).$
\end{lemma}

\begin{remark} \label{REMARK mistake in literature}
	Note that some contributions to the scientific literature [like \cite{LS:2011}, \cite{LS:2013}, \cite{BE:2019:2}] claim that the function
	\begin{align*}
		(L, R) \ni x \mapsto \int_L^{x} \Big( \mathds{1}_{(L, t]} (s) - P(t) \Big) p(s) \, \mathrm{d}s
	\end{align*}
	is differentiable when, in fact, it fails to be so in exactly one point, namely in $t$. This leads to the unfortunate consequence that we cannot assume functions in $\mathcal{F}_p$ to be differentiable, and if the random variable $X$ is discrete with an atom at the point of non-differentiability of a test function, the expectation in Lemma \ref{LEMMA density approach} makes no sense. As such, the error has no consequence for \cite{LS:2013}, since they only consider absolutely continuous random variables for which $\mathrm{disc}(X) = \emptyset$. For the general case, the restriction to test functions with $t_f \notin \mathrm{disc}(X)$ becomes necessary.
\end{remark}
\begin{remark} \label{REMARK on class F_p part 1}
	Since $f_t^{p \, \prime}(x) + \frac{p^{\prime}(x)}{p(x)} \, f^p_t(x)$ is uniformly bounded over $x \in S(p, f^p_t)$ by equation (\ref{derivative of f_t^p}), we can assume that, for each $f \in \mathcal{F}_p$, the function $f^{\prime} + \frac{p^{\prime}}{p} \, f$ is integrable with respect to any probability measure $\mathbb{P}^X$ such that $X \in S(p, f)$ $\mathbb{P}$-almost surely. Note that conditions comparable to our assumptions (C1)--(C5) are commonly stated in the context of Stein's method [see, e.g., Section 13 by \cite{CGS:2011}, Section 4 by \cite{CS:2011}, or \cite{D:2015}]. See also Remark \ref{REMARK on the regularity conditions} for further comments on the regularity conditions. It is easy to adapt the proof of Lemma \ref{LEMMA density approach} so that we can also allow for finitely many points in which the density function is zero [by changing condition (C1) accordingly]. However, for our characterization results later on we need the continuity of the functions in $\mathcal{F}_p$ on the whole interval $(L, R)$, and this we cannot get from $f_t^p$ when we allow for zeros in the function $p$.
\end{remark}
\begin{remark} \label{REMARK on class F_p part 2}
	For later use we note that if (C4) holds, any function $f \in \mathcal{F}_p$ is subject to $\lim_{x \, \searrow \, L} f(x) = 0$, since $f_t^p$ from the proof of Lemma \ref{LEMMA density approach} satisfies
	\begin{align*}
		\lim\limits_{x \, \searrow \, L} f_t^p(x)
		= \lim\limits_{x \, \searrow \, L} \frac{P(x)}{p(x)} \Big( 1 - P(t) \Big)
		= 0.
	\end{align*}
	By analogy, if (C5) holds, each function $f \in \mathcal{F}_p$ can be taken to satisfy $\lim_{x \, \nearrow \, R} f(x) = 0$.
\end{remark}
In a different form, the characterization given in Lemma \ref{LEMMA density approach} has successfully been applied for distributional approximations in the Curie-Weiss model [see \cite{CS:2011}] or the hitting times of Markov chains [see \cite{PR:2011}]. For an overview, we refer to Section 13 in \cite{CGS:2011}. In this paper, however, we use the characterization to derive another, more explicit identity that typifies distributions with density functions as above. We thereby generalize the fixed point properties of the well-known zero-bias- and equilibrium transformations but also classical identities, such as the characterization of the exponential distribution through the mean residual life function.

\section{Univariate distributions supported by the real line}
\label{SEC Real line}

Assume for now that $p : \R \to [0, \infty)$ is a probability density function supported by the whole real line.
\begin{theorem} \label{THM chara on the real line}
	Suppose that $p$ is a probability density function with $\mathrm{spt}(p) = \R$ that satisfies the conditions (C1) -- (C3). Let $X : \Omega \to \R$ be a random variable with $\mathbb{P}\big(X \in S(p)\big) = 1$, and
	\begin{align} \label{Integrability condition first theorem}
		\E \left| \frac{p^{\prime}(X)}{p(X)} \right| < \infty, \quad
		\E \left| \frac{p^{\prime}(X)}{p(X)} \, X \right| < \infty.
	\end{align}
	Then $X \sim p \mathcal{L}^1$ if, and only if, the distribution function of $X$ has the form
	\begin{align*}
		F_X(t) = \mathbb{E} \left[ \frac{p^{\prime}(X)}{p(X)} \, (t - X) \, \mathds{1}\{ X \leq t \} \right] , \quad t \in \R.
	\end{align*}
\end{theorem}
The proof is given in Appendix \ref{APP proof of real line chara}.
\begin{remark} \label{REMARK Fixed point interpretation}
	For a density function $p$ supported by the whole real line which satisfies conditions (C1) -- (C3), take the set of all distributions considered in Theorem \ref{THM chara on the real line}, that is,
	\begin{align*}
		\mathcal{P} = \left\{ \mathbb{P}^X ~ \Big| ~ \mathbb{P}\big( X \in S(p) \big) = 1,~ \E \left| \frac{p^{\prime}(X)}{p(X)} \right| < \infty, ~ \text{and} ~ \E \left| \frac{p^{\prime}(X)}{p(X)} \, X \right| < \infty \right\} .
	\end{align*}
	The previous theorem concerns properties of the mapping
	\begin{align*}
		\mathcal{T} : \mathcal{P} \to D(\R), \quad F_X \mapsto \mathcal{T}(F_X) = \left( t \mapsto \mathbb{E} \left[ \frac{p^{\prime}(X)}{p(X)} \, (t - X) \, \mathds{1}\{ X \leq t \} \right] \right),
	\end{align*}
	where $D(\R)$ is the c\`{a}dl\`{a}g-space over $\R$, and where we identified elements from $\mathcal{P}$ with their distribution function. In particular, Theorem \ref{THM chara on the real line} states that this mapping has a unique fixed point, namely $\mathbb{P}^X = p \mathcal{L}^1$. Putting further restrictions on the distribution of $X$ such that $d_p^X$ from the proof of Theorem \ref{THM chara on the real line} (see Appendix \ref{APP proof of real line chara}) is a probability density function without assuming that $F_X$ is given through our explicit formula, we have actually shown in the last calculation of that proof the existence of a distribution for some random variable $X_p$ with
	\begin{align*}
		\E \big[ f^{\prime}(X_p)\big] = \E \left[- \frac{p^{\prime}(X)}{p(X)} \, f(X)\right]
	\end{align*}
	for each $f \in \mathcal{F}_p$, and we could think of $\mathcal{T}$ as a distributional transformation. These additional restrictions [for the normal distribution they are $\E X = 0$ and $\mathrm{Var}(X) = \sigma^2$, see Example \ref{EXAMPLE Gaussian} below] scale down the class of distributions in which the characterization holds. Therefore, our point is not to cling on to distributional transformations, which makes explicit formulae more complicated [as witnessed by \cite{D:2017}, Remark 1 (d) and Remark 2], but extract whichever information we can get from the explicit formula itself.
\end{remark}

In the proof of Theorem \ref{THM chara on the real line} we have actually also shown another characterization result, but via the density function.
\begin{corollary} \label{COR chara on the real line via density function}
	Let $p$ be a probability density function with $\mathrm{spt}(p) = \R$ that satisfies conditions (C1) -- (C3). When $X : \Omega \to \R$ is a random variable with density function $\mathfrak{f}_X$, $\E \left| \tfrac{p^{\prime}(X)}{p(X)} \right| < \infty$, and $\E \left| \tfrac{p^{\prime}(X)}{p(X)} \, X \right| < \infty$, then $X \sim p \mathcal{L}^1$ if, and only if, the density function of $X$ has the form
	\begin{align*}
	\mathfrak{f}_X(t) = \mathbb{E} \left[ \frac{p^{\prime}(X)}{p(X)} \, \mathds{1}\{ X \leq t \} \right] , \quad t \in \R.
	\end{align*}
\end{corollary}

It is clear from the proof that it suffices to have the above representation for the density function of $X$ only for $\mathcal{L}^1$-almost every (a.e.) $t \in \R$ to conclude that $X \sim p \mathcal{L}^1$. This is much in line with the intuition about density functions, since they uniquely determine a probability law, but are themselves only unique $\mathcal{L}^1$-almost everywhere (a.e.). \\

To get a feeling for the results, we consider two examples. For brevity we only give the characterization via Theorem \ref{THM chara on the real line}, the result via Corollary \ref{COR chara on the real line via density function} being clear from that.
\begin{example}[\textbf{Mean-zero Gaussian distribution}] \label{EXAMPLE Gaussian}
	For $x \in \R$ let
	\begin{align*}
		p(x) = \frac{1}{\sqrt{2 \pi} \sigma}\exp\Big(-\frac{x^2}{2 \sigma^2}\Big),
	\end{align*}
	where $0 < \sigma < \infty$. The function $p$ is positive and continuously differentiable on the whole real line, so (C1) is satisfied [with $m = 0$ and $S(p) = \R$]. We have $\tfrac{p^\prime(x)}{p(x)} = - \frac{x}{\sigma^2}$, $x \in \R$. Condition (C3) follows from the existence of mean and variance of the normal distribution, and (C2) is proven using the (easily verified) identities
	\begin{align*}
		\frac{1 - P(x)}{p(x)}
		\leq \frac{\sigma^2}{x}, \quad x > 0,
		\quad \text{and} \quad
		\frac{P(x)}{p(x)}
		= \frac{1 - P(-x)}{p(-x)}
		\leq -\frac{\sigma^2}{x}, \quad x < 0.
	\end{align*}
	By Theorem \ref{THM chara on the real line}, a real-valued random variable $X$ with $\E X^2 < \infty$ follows the mean-zero Gaussian law with variance $\sigma^2$ if, and only if, the distribution function of $X$ has the form
	\begin{align*}
		F_X(t) = \mathbb{E} \left[ \frac{X}{\sigma^2} \, (X - t) \, \mathds{1}\{ X \leq t \} \right] , \quad t \in \R.
	\end{align*}
	
	In this particular example, the map $\mathcal{T}$ introduced in Remark \ref{REMARK Fixed point interpretation} is, up to a change of the domain, the zero-bias transformation discussed in the introduction. The transformation $\mathbb{P}^X \mapsto \mathbb{P}^{X^*}$ (using notation from our introduction), which coincides with our mapping $\mathcal{T}$ in terms of the law of the maps, has the normal distribution $\mathcal{N}(0, \sigma^2)$ as its unique fixed point and thus typifies this distribution within all distributions with mean zero and variance $\sigma^2$. The message of the example at hand is that our characterization result (Theorem \ref{THM chara on the real line}) has the characterization via the zero-bias distribution as a special case. It is notable that we generalize this well-known characterization in the sense that the explicit formula given above identifies the normal distribution $\mathcal{N}(0, \sigma^2)$ not only within the class of all distributions with mean zero and variance $\sigma^2$, but within the class of all distributions with $\E X^2 < \infty$. However, if $\E X \neq 0$ or $\mathrm{Var}(X) \neq \sigma^2$, the formula for $F_{X^*}$ may no longer be a distribution function, and $\mathcal{T}$ is to be understood as an extension of the operator that maps $\mathbb{P}^X \mapsto \mathbb{P}^{X^*}$ onto the larger domain
	\begin{align*}
	\mathcal{P}
	= \Big\{ \mathbb{P}^X ~ \Big| ~ \E X^2 < \infty \Big\}
	\supsetneq \Big\{ \mathbb{P}^X ~ \Big| ~ \E X = 0 ~ \text{and} ~ \mathrm{Var}(X) = \sigma^2 \Big\}.
	\end{align*}
	
	The conditions (C1) -- (C3) also hold for the normal distribution with location parameter included. We simply chose the setting above to illustrate the connection to the zero-bias distribution.
\end{example}

\begin{example}[\textbf{Laplace distribution}] \label{EXAMPLE Laplace}
	For a location parameter $\mu \in \R$ and a scale parameter $\sigma > 0$ consider the density of the corresponding Laplace distribution,
	\begin{align*}
	p(x) = \frac{1}{2 \sigma} \exp\Big( - \frac{|x - \mu|}{\sigma} \Big), \quad x \in \R.
	\end{align*}
	Condition (C1) is satisfied with $m = 1$, $y_1 = \mu$, and $S(p) = (- \infty, \mu) \cup (\mu, \infty)$. We have
	\begin{align*}
		\frac{p^\prime(x)}{p(x)} = \frac{\mathrm{sign}(\mu - x)}{\sigma}, \quad x \neq \mu.
	\end{align*}
	To verify (C2), use that the distribution function of the Laplace distribution can be given explicitly to obtain $\sup_{x \, \in \, S(p)} \kappa_p(x) \leq 1 < \infty$. Condition (C3) follows from a simple calculation. Consequently, Theorem \ref{THM chara on the real line} holds, and the characterization for the Laplace distribution reads as follows. A real-valued random variable $X$ with distribution function $F_X$, which satisfies $\mathbb{P}\big( X \neq \mu \big) = 1$ and $\E |X| < \infty$, has the Laplace distribution with parameters $\mu$ and $\sigma$ if, and only if,
	\begin{align*}
		F_X(t) = \mathbb{E} \left[ \frac{\mathrm{sign}(\mu - X)}{\sigma} \, (t - X) \, \mathds{1}\{ X \leq t \} \right] , \quad t \in \R.
	\end{align*}
\end{example}
\vspace{5mm}
In the context of probability distributions on the real line, we have also checked the conditions (C1) -- (C3) for the Cauchy- and Gumbel distribution, showing that we do not need any moment assumptions to prove (C3), and that the characterizations include more complicated distributions which are important in applications. We will give more examples later on.

\section{Univariate distributions with semi-bounded support}
\label{SEC semi-bounded support}
In this section, we seek to provide characterization results similar to those in the previous section, but for probability distributions with semi-bounded support. We have chosen in Section \ref{SEC Real line} to first prove the characterization via the distribution function since this is the 'conventional' way, or at least, say, the way the special case of the zero-bias transformation is known. From a logical perspective it is more convenient to first establish the result via the density function as in Corollary \ref{COR chara on the real line via density function}, and then to derive the corresponding distribution function. We first discuss the case when $p$ is a density functions whose support is bounded from below. Namely, we let $p : \R \to [0, \infty)$ be a probability density function with $\mathrm{spt}(p) = [L, \infty)$, $L > - \infty$. The most important case is $L = 0$, that is, density functions supported by the positive half line.
\begin{theorem} \label{THM chara semi-bounded support via density function}
	Let $p$ be a probability density function with $\mathrm{spt}(p) = [L, \infty)$ that satisfies the conditions (C1) -- (C4). If $X : \Omega \to (L, \infty)$ is a random variable with density function $\mathfrak{f}_X$, $\E \left| \tfrac{p^{\prime}(X)}{p(X)} \right| < \infty$, and $\E \left| \tfrac{p^{\prime}(X)}{p(X)} \, X \right| < \infty$, then $X \sim p \mathcal{L}^1$ if, and only if,
	\begin{align*}
		\mathfrak{f}_X(t)
		= \mathbb{E} \left[ - \frac{p^{\prime}(X)}{p(X)} \, \mathds{1}\{ X > t \} \right] , \quad t > L.
	\end{align*}
\end{theorem}
The proof of this theorem consists of arguments and calculations that are very similar to those in the proof of Theorem \ref{THM chara on the real line}, and we refrain from giving the details. Instead we give some insight on the special case of density functions on the positive axis.

\begin{remark} \label{REMARK integr. condition in sufficiency part of density char. semi-bounded sup.}
	The integrability condition on $X$ can be weakened in cases where the density function $p$ is positive and continuously differentiable, as well as supported by the positive axis, that means especially, $m = 0$ and $S(p) = (0, \infty)$. In this case, the calculation in the sufficiency part of the proof of Theorem \ref{THM chara semi-bounded support via density function} reduces to
	\begin{align*}
		\E\big[ f^{\prime}(X) \big]
		&= \int_{0}^{t_f} f^{\prime}(s) \, \E\left[ - \frac{p^{\prime}(X)}{p(X)} \, \mathds{1}\{X > s\} \right] \mathrm{d}s + \int_{t_f}^{\infty} f^{\prime}(s) \, \E\left[ - \frac{p^{\prime}(X)}{p(X)} \, \mathds{1}\{X > s\} \right] \mathrm{d}s \\
		&= \E\left[ - \frac{p^{\prime}(X)}{p(X)} \int_{0}^{X} f^{\prime}(s) \, \mathrm{d}s \, \mathds{1}\{ X \leq t_f \} \right] + \E\left[ - \frac{p^{\prime}(X)}{p(X)} \int_{0}^{t_f} f^{\prime}(s) \, \mathrm{d}s \, \mathds{1}\{ X > t_f \} \right] \\
		&~~~ + \E\left[ - \frac{p^{\prime}(X)}{p(X)} \int_{t_f}^{X} f^{\prime}(s) \, \mathrm{d}s \, \mathds{1}\{ X > t_f \} \right] \\
		&= \E\left[ -\frac{p^{\prime}(X)}{p(X)} \, f(X) \right],
	\end{align*}
	and it suffices for the use of Fubini's theorem to know that $\E \left| \tfrac{p^{\prime}(X)}{p(X)} \, X \right| < \infty$ (note that this condition on $X$ is also enough to guarantee that the expectation which defines $d_p^X$ exists $\mathcal{L}^1$-a.e., see Appendix \ref{APP proof semi-bounded supp chara}). Consequently, it suffices to claim $\int_0^\infty x |p^\prime(x)| \, \mathrm{d}x < \infty$ instead of (C3), and still have the sufficiency part of the theorem be consistent in itself. What is more, this last condition yields
	\begin{align*}
		\int_0^\infty \int_t^\infty \big| p^\prime(x) \big| \, \mathrm{d}x \, \mathrm{d}t
		= \int_0^\infty \big| p^\prime(x) \big| \int_0^x \mathrm{d}t \, \mathrm{d}x
		= \int_0^{\infty} x \big|p^{\prime}(x)\big| \, \mathrm{d}x
		< \infty,
	\end{align*}
	and thus $\int_t^{\infty} \big|p^{\prime}(x)\big| \, \mathrm{d}x < \infty$ for $\mathcal{L}^1$-a.e. $t > 0$. This suffices to derive the necessity part of Theorem \ref{THM chara semi-bounded support via density function} with equality for $\mathcal{L}^1$-a.e. $t > 0$. Putting together these thoughts, we obtain the following special case.
\end{remark}

\begin{corollary}[\textbf{Densities supported by the positive axis}] \label{COR chara supported by pos. half axis, via density function}
	Assume that $p$ is a probability density function with $\mathrm{spt}(p) = [0, \infty)$ that is positive and continuously differentiable on $(0, \infty)$, and satisfies (C2) and (C4). Moreover, assume that $\int_0^{\infty} x \big|p^{\prime}(x)\big| \, \mathrm{d}x < \infty$. Let $X$ be a positive random variable with density function $\mathfrak{f}_X$, and $\E \left| \tfrac{p^{\prime}(X)}{p(X)} \, X \right| < \infty$. Then $X \sim p \mathcal{L}^1$ if, and only if, we have for $\mathcal{L}^1$-a.e. $t > 0$ that
	\begin{align*}
		\mathfrak{f}_X (t) = \mathbb{E} \left[- \frac{p^{\prime}(X)}{p(X)} \, \mathds{1}\{ X > t \} \right] .
	\end{align*}
\end{corollary}

Up next, we use Theorem \ref{THM chara semi-bounded support via density function} to derive a characterization result for the distribution function.
\begin{theorem} \label{THM chara distr. function, sup. bounded from below}
	Assume that $p$ is a probability density function with $\mathrm{spt}(p) = [L, \infty)$ satisfying the conditions (C1) -- (C4). Let $X : \Omega \to (L, \infty)$ be a random variable with $\mathbb{P}\big( X \in S(p) \big) = 1$, $\E \left| \frac{p^{\prime}(X)}{p(X)} \right| < \infty$, and $\E \left| \frac{p^{\prime}(X)}{p(X)} \, X \right| < \infty$. Then $X \sim p \mathcal{L}^1$ if, and only if,
	\begin{align*}
		F_X (t) = \E \left[ - \frac{p^{\prime}(X)}{p(X)} \Big( \min\{X, t\} - L \Big) \right], \quad t > L.
	\end{align*}
\end{theorem}
\vspace{5mm}
The proof is given in Appendix \ref{APP proof semi-bounded supp chara}. Note that the results on the distribution function are somewhat richer than the characterizations via the density function, for the latter only identify the underlying distribution within a subset of absolutely continuous probability distributions for which a density function exists. The characterization via the distribution function does not need this restriction to absolutely continuous distributions, but only that $X$ has no atoms in $(L, \infty) \setminus S(p) = \{ y_1, \dots, y_m \}$.

\begin{remark} \label{REMARK alternative condition (C3) distr. chara semibounded case}
	In the case where $L = 0$, and $p$ is continuously differentiable and positive on $(0, \infty)$, the proof of Theorem \ref{THM chara distr. function, sup. bounded from below} remains true if we replace (C3) with $\int_0^{\infty} x |p^{\prime}(x)| \, \mathrm{d}x < \infty$, and if we further drop the first integrability condition on $X$ and only require $\E \left| \tfrac{p^\prime(X)}{p(X)} \, X \right| < \infty$.
\end{remark}

We obtain the following special case of the characterization.
\begin{corollary}[\textbf{Densities supported by the positive axis}] \label{COR chara supported by pos. half axis, via distribution function}
	Assume that $p$ is a probability density function with $\mathrm{spt}(p) = [0, \infty)$ that is positive and continuously differentiable on $(0, \infty)$, and satisfies (C2) and (C4). Moreover, assume that $\int_0^{\infty} x \big|p^{\prime}(x)\big| \, \mathrm{d}x < \infty$. Let $X$ be a positive random variable with $\E \left| \tfrac{p^{\prime}(X)}{p(X)} \, X \right| < \infty$. Then $X \sim p \mathcal{L}^1$ if, and only if,
	\begin{align*}
		F_X(t) = \E \left[ - \frac{p^{\prime}(X)}{p(X)} \, \min\{X, t\} \right], \quad t > 0.
	\end{align*}
\end{corollary}

Now follows the major source of examples we give in this work. We omit the explicit proofs of the regularity conditions for they consist of (sometimes) tedious calculations which provide no insight on the characterizations. Instead, we give the following remark on how the conditions are to be verified, and on their necessity in general.
\begin{remark} \label{REMARK on the regularity conditions}
	The regularity condition (C1) is easily understood and checked for a given density function. Note that the weaker assumption of absolute continuity of $p$, which is mostly used in the context of Stein's method, entails similar problems as described in Remark \ref{REMARK mistake in literature}: If $p$ is merely assumed to be absolutely continuous, then in order to handle random variables $X$ with discrete parts (e.g. in Lemma \ref{LEMMA density approach}) we would still have to identify the points of non-differentiability of $p$ in order to make sense of the term $p^\prime(X)$. This would return us to considering a set like $S(p)$ which, technically, brings us to the setting we consider already.
	
	Condition (C3) involves a direct calculation which can often be simplified if one has knowledge of the existence of moments of the distribution at hand. From the proofs of our characterizations it is apparent that (C1) and (C3) [as well as the integrability conditions on $X$ which are in line with (C3)] are necessary to use Fubini's theorem and the fundamental theorem of calculus. As such, we do not see any truly instrumental weaker alternative conditions (apart from the special case discussed in Remarks \ref{REMARK integr. condition in sufficiency part of density char. semi-bounded sup.} and \ref{REMARK alternative condition (C3) distr. chara semibounded case}) which still rigorously allow for all calculations.
	
	Both conditions (C4) and (C5) are trivially satisfied when the respective limit of the density function is positive, and if that is not the case, L'Hospital's rule gives a reliable handle for it. With regard to these two conditions, we refer to Proposition 3.7 of \cite{D:2015} who discusses them in much detail and provides easy-to-check criteria. Moreover, this specific result from \cite{D:2015} indicates strongly that the two conditions are not restrictive in practice.
	
	To prove condition (C2), it is helpful to realize, in the case when $p$ is continuously differentiable, that $\kappa_p$ is continuous. Thus, it suffices to check that $\limsup_{x \, \searrow \, L} \kappa_p(x) < \infty$ and $\limsup_{x \, \nearrow \, R} \kappa_p(x) < \infty$ for (C2) to hold. Regularity conditions (C2) and (C4)/(C5) guarantee certain beneficial properties of the test functions from $\mathcal{F}_p$. For one, they guarantee that, for $f \in \mathcal{F}_p$, $\lim_{x \, \searrow \, L} f(x) = 0$ (or $\lim_{x \, \nearrow \, R} f(x) = 0$), see Remark \ref{REMARK on class F_p part 2}, which we need to truly get rid of the test functions in our calculations (as in Appendix \ref{APP proof bounded supp chara}). Condition (C2) is stated so that functions $f \in \mathcal{F}_p$ have uniformly bounded derivative. We use this fact in our proofs (e.g., the last calculation in Appendix \ref{APP proof of real line chara}) to apply the fundamental theorem of calculus on $f^\prime$ and to justify the use of Fubini's theorem. For both arguments the boundedness of $f^\prime$ is not a necessary condition, but we have not found any alternative assumption for (C2) which allows for a sound and rigorous derivation of all results.
	
	Later on, we give an example for a distribution which fails the respective version of condition (C3) that ought to hold in order for that distribution to be included in our characterization results. For a (rather artificial) density functions which violates (C4), see Example 3.6 of \cite{D:2015}.
\end{remark}
\vspace{3mm}
With these tools at hand, the regularity conditions for all examples below can be proven. We use Corollary \ref{COR chara supported by pos. half axis, via distribution function} in each case, except for the L\'{e}vy distribution. The characterizations via the density functions are not stated explicitly to save space.
\begin{example}[\textbf{Gamma distribution}] \label{EXAMPLE Gamma}
	Assume that
	\begin{align*}
	p(x) = \frac{\lambda^{-k}}{\Gamma(k)} \, x^{k - 1} \, \exp\big( - \lambda^{-1} x \big), \quad x > 0,
	\end{align*}
	is the density function of the Gamma distribution with shape parameter $k > 0$ and scale parameter $\lambda > 0$. If $X$ is a positive random variable with $\E X < \infty$, then $X$ follows the Gamma law with parameters $k$ and $\lambda$ if, and only if, the distribution function of $X$ has the form
	\begin{align*}
		F_X(t)
		= \E \left[ \left( - \frac{k - 1}{X} + \frac{1}{\lambda} \right) \min\{ X, t \} \right], \quad t > 0.
	\end{align*}
	Note that this result has been proven explicitly, and with a similar line of proof as our general results above, by \cite{BE:2019:2}.
\end{example}
\begin{example}[\textbf{Exponential Distribution}] \label{EXAMPLE Exponential}
	Denote the density of the exponential distribution with rate parameter $\lambda > 0$ by $p(x) = \lambda e^{-\lambda x}$, $x > 0$.	This is an easy special case of the previous example, namely the Gamma distribution with shape parameter $k = 1$ and scale parameter $1 / \lambda$. Let $X$ be a positive random variable with $\E X < \infty$. Then $X$ has the exponential distribution with parameter $\lambda$ if, and only if,
	\begin{align*}
		F_X(t)
		= \lambda \, \E \Big[ \min\{ X, t \} \Big], \quad t > 0.
	\end{align*}

	This identity is [see \cite{BH:2000}] equivalent to the well-known characterization of exponentiality via the mean residual life function, which states that a positive random variable $X$ with $\E X < \infty$ follows an exponential law if, and only if, $\E [ X - t \,|\, X > t ] = \E[X]$, $t > 0$. For yet another observation, assume that $X$ is a positive random variable with $\E X = \lambda^{-1}$. With
	\begin{align*}
		d_p^X(t)
		= \E \left[ - \frac{p^\prime(X)}{p(X)} \, \mathds{1}\{ X > t \} \right]
		= \lambda \, \mathbb{P}\big( X > t \big), \quad t > 0,
	\end{align*}
	as in the proofs of our results, we have $d_p^X \geq 0$ and
	\begin{align*}
		\int_0^{\infty} d_p^X(t) \, \mathrm{d}t
		= \lambda \int_0^{\infty} \mathbb{P}\big(X > t\big) \, \mathrm{d}t
		= \lambda \, \E X
		= 1.
	\end{align*}
	If $X^{e}$ is a random variable with density function $d_p^X$, the proof of Theorem \ref{THM chara semi-bounded support via density function} (see Remark \ref{REMARK alternative condition (C3) distr. chara semibounded case}) shows that $\E [ f^{\prime}(X^{e}) ] = \lambda \, \E[ f(X) ]$ for each $f \in \mathcal{F}_p$. Up to a change in the class of test functions, this is the defining equation of the equilibrium distribution with respect to $X$. Lemma \ref{LEMMA density approach} implies that when restricting to $\E X = \lambda^{-1}$, the exponential distribution with parameter $\lambda$ is the unique fixed point of the equilibrium transformation $\mathbb{P}^X \mapsto \mathbb{P}^{X^{e}}$. This fact is used for approximation arguments with Stein's method [see \cite{PR:2011}, who introduced the equilibrium distribution, as well as Chapter 13.4 by \cite{CGS:2011} and Section 5 by \cite{R:2011}]. As in the case of the zero-bias transformation, we have generalized this characterization in the sense that the explicit formula of the equilibrium distribution uniquely identifies the exponential distribution with parameter $\lambda$ within the class of all distributions $\mathbb{P}^X$ with $\E X < \infty$.
\end{example}

\begin{example}[\textbf{Inverse Gaussian distribution}] \label{Example inverse Gaussian}
	Denote the inverse Gaussian density by
	\begin{align*}
		p(x) = \sqrt{\frac{\lambda}{2 \pi}} \, x^{- 3/2} \exp\left( - \frac{\lambda (x - \mu)^2}{2 \mu^2 x} \right), \quad x > 0,
	\end{align*}
	where $\mu, \lambda > 0$. If $X$ is a positive random variable with $\E\big[ X + X^{-1} \big] < \infty$, then $X$ follows the inverse Gaussian law with parameters $\mu$ and $\lambda$ if, and only if,
	\begin{align*}
		F_X(t)
		= \E \left[ \left( - \frac{\lambda}{2 X^2} + \frac{3}{2 X} + \frac{\lambda}{2 \mu^2} \right) \min\{ X, t \} \right], \quad t > 0.
		\end{align*}
\end{example}
\vspace{5mm}
Now we handle distributions that are of interest for applications. The Weibull distribution is applied in hydrology and wind speed analysis, see \cite{S:1987} and \cite{CCDO:2014}, the Burr distribution is commonly taken as a model for household income, see \cite{SM:1976}, and the Rice distribution appears in signal processing to describe how cancellation phenomena affect radio signals [cf. Chapter 13 of \cite{PS:2008}]. The last example we give is the L\'{e}vy distribution which is used to model the length of paths that are followed by photons after reflection from a turbid media, see Section 3 of \cite{R:2008}. Here we provide insight on the handling of an additional location parameter which is often added to probability distributions.

\begin{example}[\textbf{Weibull distribution}] \label{EXAMPLE Weibull}
	For $k, \lambda > 0$ let
	\begin{align*}
		p(x) = \frac{k}{\lambda^k} \, x^{k - 1} \exp\left( - \Big( \frac{x}{\lambda} \Big)^k \right), \quad x > 0,
	\end{align*}
	be the density function of the Weibull distribution in its usual parametrization. Let $X$ be any positive random variable with $\E X^k < \infty$. Then $X$ has the Weibull distribution with parameters $k$ and $\lambda$ if, and only if,
	\begin{align*}
		F_X(t)
		= \E \left[ \left( \frac{k \, X^{k - 1}}{\lambda^k} - \frac{k - 1}{X} \right) \min\{ X, t \} \right], \quad t > 0.
	\end{align*}
\end{example}

\begin{example}[\textbf{Burr distribution}] \label{EXAMPLE Burr}
	The Burr Type XII distribution with parameters $c, k > 0$ and scale parameter $\sigma > 0$ has density function
	\begin{align*}
		p(x)
		= \frac{c \, k}{\sigma} \, \Big( \frac{x}{\sigma} \Big)^{c - 1} \left( 1 + \Big( \frac{x}{\sigma} \Big)^c \right)^{-k - 1}, \quad x > 0.
	\end{align*}
	A positive random variable $X$ has the Burr distribution with parameters $c, k, \sigma > 0$ if, and only if, the distribution function of $X$ has the form
	\begin{align*}
		F_X(t)
		= \E \left[ \left( c \, (k + 1) \, \frac{X^{c - 1}}{\sigma^c + X^c} - \frac{c - 1}{X} \right) \min\{ X, t \} \right], \quad t > 0.
	\end{align*}
	Particularly interesting about this example is that, even though the Burr distribution is substantially more complicated than many of our other examples, no moment condition is needed for the characterization to hold, since
	\begin{align*}
		\E \left| \frac{p^\prime(X)}{p(X)} \, X \right|
		\leq |c - 1| + c \, (k + 1) \, \E \left[ \frac{X^c}{\sigma^c + X^c} \right]
		\leq |c - 1| + c \, (k + 1)
		< \infty.
	\end{align*}
	This implies that the characterization is universal in the sense that it identifies the Burr distribution within the set of all probability distributions on the positive axis.
\end{example}

\begin{example}[\textbf{Rice distribution}] \label{EXAMPLE Rice}
	For parameters $k, \varrho > 0$ the density function of the Rice distribution is given by
	\begin{align*}
		p(x)
		= \frac{2 \, (k + 1) \, x}{\varrho} \, \exp\left( - k - \frac{(k + 1) \, x^2}{\varrho} \right) I_0\left( 2 \, \sqrt{\frac{k \, (k + 1)}{\varrho}} \, x \right), \quad x > 0,
	\end{align*}
	where $I_\alpha$ denotes the modified Bessel function of first kind of order $\alpha \in \Z$. We chose the parametrization for $p$ that is mostly used in signal processing and is easily found under the keyword of Rician fading. Let $X$ be a positive random variable with $\E X^2 < \infty$. Then $X$ has the Rice distribution with parameters $k$ and $\varrho$ if, and only if
	\begin{align*}
		F_X(t)				
		= \E\left[ \left( - \frac{1}{X} + \frac{2 \, (k + 1) \, X}{\varrho} - 2 \, \sqrt{\frac{k \, (k + 1)}{\varrho}} \cdot \frac{I_1\left( 2 \, \sqrt{\tfrac{k \, (k + 1)}{\varrho}} \, X \right)}{I_0\left( 2 \, \sqrt{\tfrac{k \, (k + 1)}{\varrho}} \, X \right)} \right) \min\{ X, t \} \right],
	\end{align*}
	for $t > 0$. Note that despite the complexity of the term $\tfrac{p^\prime(x)}{p(x)}$, the integrability conditions is $\E X^2 < \infty$, since the quotient of the Bessel functions cancels via $\tfrac{I_1(y)}{I_0(y)} \leq 1$, $y > 0$.
\end{example}

\begin{example}[\textbf{L\'{e}vy Distribution}] \label{EXAMPLE Levy}
	Take $\mu \in \R$ and $\sigma > 0$. Let
	\begin{align*}
		p(x) = \sqrt{\frac{\sigma}{2 \pi}} \, \big( x - \mu \big)^{- 3/2} \exp\left( - \frac{\sigma}{2 (x - \mu)} \right), \quad x > \mu,
	\end{align*}
	denote the density function of the L\'{e}vy distribution with location parameter $\mu$ and scale parameter $\sigma$. Let $X$ be a random variable which takes values in $(\mu, \infty)$ almost surely such that $\E [ ( X - \mu )^{-1} ] < \infty$ and $\E [ ( X - \mu )^{-2} ] < \infty$. Then $X$ has the L\'{e}vy distribution with parameters $\mu$ and $\sigma$ if, and only if, the distribution function of $X$ has the form
	\begin{align*}
		F_X(t)
		= \frac{1}{2} \, \E \left[ \left( \frac{3}{X - \mu} - \frac{\sigma}{(X - \mu)^2} \right) \Big( \min\{ X, t \} - \mu \Big) \right], \quad t > \mu.
	\end{align*}
\end{example}
\vspace{5mm}
The following example is one which fails the regularity condition (C3). Recall that for distributions which are not supported by the positive axis, we need (C3) fully, that is, we cannot apply Remarks \ref{REMARK integr. condition in sufficiency part of density char. semi-bounded sup.} or \ref{REMARK alternative condition (C3) distr. chara semibounded case}.
\begin{example}[\textbf{Shifted Gamma distribution}] \label{EXAMPLE shifted gamma distribution}
	Assume that
	\begin{align*}
		p(x) = \frac{\lambda^{-k}}{\Gamma(k)} \, (x - \mu)^{k - 1} \, \exp\big( - \lambda^{-1} (x - \mu) \big), \quad x > \mu,
	\end{align*}
	is the density function of the shifted Gamma distribution with shape parameter $k > 0$, scale parameter $\lambda > 0$, and location parameter $\mu \in \R \setminus \{0\}$. We have
	\begin{align*}
		\frac{p^\prime(x)}{p(x)}
		= \frac{k - 1}{x - \mu} - \frac{1}{\lambda}, \quad x > \mu.
	\end{align*}
	Since $\mu \neq 0$, in order to establish our characterization result, we have to verify the conditions from Theorem \ref{THM chara distr. function, sup. bounded from below} which includes (C3). However, for $k < 1$ we have
	\begin{align*}
		\int_\mu^\infty \big| p^\prime(x) \big| \, \mathrm{d}x
		&\geq \int_\mu^\infty \frac{|k - 1|}{x - \mu} \,  p(x) \, \mathrm{d}x - \frac{1}{\lambda} \int_\mu^\infty p(x) \, \mathrm{d}x \\
		&= \frac{1 - k}{\lambda \, \Gamma(k)} \int_0^\infty z^{k - 2} \, e^{- z} \, \mathrm{d}z - \frac{1}{\lambda} \\
		&= \infty.
	\end{align*}
\end{example}

Next, we discuss the characterizations for probability distributions supported by the positive axis in the case of exponential families. More specifically, we focus on continuously differentiable density functions. Quite a few of the examples we already gave can be written as an exponential family, but we do not reconsider them and instead give a new example at the end of this part. Of course the arguments below could also be used to treat exponential families over the real line, using Theorem \ref{THM chara on the real line}. In detail, we let $\Theta \subset \R^d$ be non-empty, and consider an exponential family (over the positive axis) in the natural parametrization given through
\begin{align*}
	p_\vartheta(x)
	= c(\vartheta) \, h(x) \exp\Big( \vartheta^\top T(x) \Big), \quad x > 0, \, \vartheta \in \Theta,
\end{align*}
where $T = (T_1, \dots, T_d)^\top : (0, \infty) \to \R^d$ and $h : (0, \infty) \to [0, \infty)$ are (Borel-) measurable functions, $\vartheta^\top$ is the transpose of a column vector $\vartheta$, and 
\begin{align*}
	c(\vartheta)
	= \left( \int_0^\infty h(x) \exp\Big( \vartheta^\top T(x) \Big) \mathrm{d}x \right)^{-1} .
\end{align*}
We choose $\Theta$ such that $0 < c(\vartheta) < \infty$ for each $\vartheta \in \Theta$. The exponential family is assumed to be strictly $d$-parametric, that is, we take the functions $1, T_1, \dots, T_d$ to be linearly independent on the complement of every null set. The definition of exponential families, and insights on their properties, are provided by virtually any classical textbook on mathematical statistics.

We try to get an idea on how the conditions (C1) -- (C4) can be handled for exponential families. Condition (C2) remains a little cryptic, meaning that it depends on the given example how it can be proven, and, at this point, we cannot give any improvement to what we discussed in Remark \ref{REMARK on the regularity conditions} concerning that condition.
\begin{itemize}
	\item[\textbf{(C1)}]
	{
		Assume that $T$ and $h$ are continuously differentiable, and that $h$ is positive. Trivially, these assumptions cover (C1) for they assure that for each $\vartheta \in \Theta$, $p_\vartheta$ is continuously differentiable and positive on $(0, \infty)$. For $x > 0$ we have
		\begin{align*}
			\frac{p_\vartheta^\prime(x)}{p_\vartheta(x)}
			= \vartheta^\top T^\prime(x) + \frac{h^\prime(x)}{h(x)},
		\end{align*}
		where $T^\prime(x) = \big( T_1^\prime(x), \dots, T_d^\prime(x) \big)^\top$.
	}
	\item[\textbf{(C3)}]
	{
		Using the weaker subsidy for (C3) given in the Remarks \ref{REMARK integr. condition in sufficiency part of density char. semi-bounded sup.} and \ref{REMARK alternative condition (C3) distr. chara semibounded case}, a sufficient condition for (C3) is derived as follows. Let $\vartheta \in \Theta$, and take $Z \sim p_\vartheta \mathcal{L}^1$. Then
		\begin{align*}
			\int_0^\infty x \big| p_\vartheta^\prime(x) \big| \, \mathrm{d}x
			\leq \big\lVert \vartheta \big\rVert \, \E \Big[ \big\lVert T^\prime(Z) \big\rVert Z \Big] + \E \left[ \left| \frac{h^\prime(Z)}{h(Z)} \right| Z \right].
		\end{align*}
		Therefore, it suffices to know that
		\begin{align} \label{integrability condition (C3) for expon. family}
			\E \left| \frac{h^\prime(Z)}{h(Z)} \, Z \right| < \infty
			\quad \text{and} \quad
			\E\Big[ \big| T_j^\prime(Z) \big| \, Z \Big] < \infty, \quad j = 1, \dots, d.
		\end{align}
		Since $T$ often consists of monomials $x^{k}$, $k \in \Z$, or of some logarithmic term $\log (x)$, (\ref{integrability condition (C3) for expon. family}) frequently reduces to a moment constraint which is satisfied if the expectation of $T(Z)$ exists.
	}
	\item[\textbf{(C4)}]
	{
		Note that $P_\vartheta$, the distribution function corresponding to $p_\vartheta$, satisfies $\lim_{x \, \searrow \, 0} P_\vartheta(x) = 0$, so if $\lim_{x \, \searrow \, 0} p_\vartheta(x) > 0$ , (C4) is trivially satisfied. If $\lim_{x \, \searrow \, 0} p_\vartheta(x) = 0$, a sufficient condition for (C4) is that
		\begin{align*}
			\lim_{x \, \searrow \, 0} \left( \vartheta^\top T^\prime(x) + \frac{h^\prime(x)}{h(x)} \right)
			= \infty.
		\end{align*}
	}
\end{itemize}
We now give the characterization result that follows from Corollary \ref{COR chara supported by pos. half axis, via distribution function}. Corollary \ref{COR chara supported by pos. half axis, via density function} yields a similar result via the density function, but we will not restate it explicitly.
\begin{corollary} \label{COR expon. families chara. pos. axis, via distribution function}
	Let $\big\{ p_\vartheta \, | \, \vartheta \in \Theta \big\}$ be an exponential family as above. Assume that each $p_\vartheta$ is continuously differentiable and positive, and satisfies (C2) -- (C4). Let $X$ be a positive random variable with
	\begin{align*}
		\E \left[ \left( \big\lVert T^\prime(X) \big\rVert + \left|\frac{h^\prime(X)}{h(X)}\right| \right) X \right] < \infty.
	\end{align*}
	Then $X \sim p_\vartheta \mathcal{L}^1$ if, and only if, the distribution function of $X$ has the form
	\begin{align*}
		F_X(t)
		= \E \left[ - \left( \vartheta^\top T^\prime(X) + \frac{h^\prime(X)}{h(X)} \right) \min\{ X, t \} \right], \quad t > 0.
	\end{align*}
\end{corollary}
\vspace{1mm}

\begin{example}[\textbf{Log-normal distribution}] \label{EXAMPLE log-normal}
	For parameters $\mu \in \R$ and $\sigma > 0$ consider the density function of the log-normal distribution
	\begin{align*}
	p(x)
	&= \frac{1}{x \, \sqrt{2 \pi} \, \sigma} \exp\left( - \frac{\big( \log(x) - \mu \big)^2}{2 \sigma^2} \right) \\
	&= \sqrt{- 2 \vartheta_2} \exp\Big( \frac{\vartheta_1^2}{4 \vartheta_2} \Big) \frac{1}{\sqrt{2 \pi} \, x} \exp\Big( \vartheta_1 \, \log(x) + \vartheta_2 \, \log^2(x) \Big), \quad x > 0,
	\end{align*}
	where $\vartheta = \big( \vartheta_1, \, \vartheta_2 \big)^\top = \big( \tfrac{\mu}{\sigma^2}, \, - \tfrac{1}{2 \sigma^2} \big)^\top$. In the last representation we see that the class of log-normal distributions forms an exponential family with parameter space $\Theta = \R \times (- \infty, 0)$, as well as $h(x) = \tfrac{1}{\sqrt{2 \pi} \, x}$, $T(x) = \big( \log(x), \, \log^2(x) \big)^\top$, and $c(\vartheta) = \sqrt{- 2 \vartheta_2} \exp\big( \tfrac{\vartheta_1^2}{4 \vartheta_2} \big)$, where $c(\vartheta) \in (0, \infty)$ for every $\vartheta \in \Theta$. In this whole example we suppress the index $\vartheta$ for $p$. All of the following arguments are valid for any fixed (but arbitrary) $\vartheta \in \Theta$.
	
	The density function $p$ is continuously differentiable since $h$ and $T$ are such, and it is positive as $h$ is so. For $x > 0$ we have
	\begin{align*}
		\frac{p^\prime(x)}{p(x)}
		= \vartheta^\top T^\prime(x) + \frac{h^\prime(x)}{h(x)}
		= \frac{\big( \mu - \sigma^2 \big) - \log(x)}{\sigma^2 \, x} .
	\end{align*}
	For the log-normal density function we have $\lim_{x \, \searrow \, 0} p(x) = 0$, as well as
	\begin{align*}
		\lim_{x \, \searrow \, 0} \left( \vartheta^\top T^\prime(x) + \frac{h^\prime(x)}{h(x)} \right)
		= \lim_{x \, \searrow \, 0} \frac{\big( \mu - \sigma^2 \big) - \log(x)}{\sigma^2 \, x}
		= \infty ,
	\end{align*}
	and the discussion of (C4) yields that this condition holds. In order to establish (C3), let $Z \sim p \mathcal{L}^1$. Then $\log(Z)$ is Gaussian with mean $\mu$ and variance $\sigma^2$, and the expectation of $\log(Z)$ exists, that is, $\E | \log(Z) | < \infty$. Therefore, we have
	\begin{align*}
		\E \left| \frac{h^\prime(Z)}{h(Z)} \, Z \right| = 1 < \infty, \quad
		\E \big| T_1^\prime(Z) \, Z \big| = 1 < \infty, \quad \text{and} \quad
		\E \big| T_2^\prime(Z) \, Z \big|
		= 2 \, \E \big| \log(Z) \big|
		< \infty,
	\end{align*}
	which suffices for (C3) by the discussions above. The proof of (C2) is a bit tedious and follows Remark \ref{REMARK on the regularity conditions}. As it provides no insight on that regularity condition, we omit it here. The characterization result for the log-normal distribution as given in the Corollary \ref{COR expon. families chara. pos. axis, via distribution function} is as follows. If $X$ is a positive random variable with $\E | \log(X) | < \infty$, then $X$ follows the log-normal law with parameters $\mu \in \R$ and $\sigma > 0$ if, and only if, the distribution function of $X$ has the form
	\begin{align*}
		F_X(t)
		= \E \left[ - \frac{( \mu - \sigma^2 ) - \log(X)}{\sigma^2X} \, \min\{ X, t \} \right], \quad t > 0.
	\end{align*}
\end{example}
\vspace{3mm}

Finally, we state the characterization result for a probability density function $p$ with support bounded from above, $\mathrm{spt}(p) = (- \infty, R]$, $R < \infty$. We omit the proof since it is a collage of earlier proofs, and only state the result. A characterization via the density function also holds, and its form is immediately conceivable from the result we state below. Similar observations concerning integrability conditions carry over from the case of density function with support bounded from below.
\begin{corollary} \label{COR chara support bounded from above via distribution function}
	Let $p$ be a probability density function with $\mathrm{spt}(p) = (- \infty, R]$ that satisfies the conditions (C1) -- (C3) and (C5). Take $X : \Omega \to (- \infty, R)$ to be a random variable with $\mathbb{P}\big( X \in S(p) \big) = 1$, $\E \left| \tfrac{p^{\prime}(X)}{p(X)} \right| < \infty$, and $\E \left| \tfrac{p^{\prime}(X)}{p(X)} \, X \right| < \infty$. Then $X \sim p \mathcal{L}^1$ if, and only if,
	\begin{align*}
		1 - F_X(t) = \E \left[ \frac{p^{\prime}(X)}{p(X)} \Big( R - \max\{ X, t \} \Big) \right], \quad t < R.
	\end{align*}	
\end{corollary}

\section{Univariate distributions with bounded support}
\label{SEC bounded support}

For the sake of completeness, we study probability density functions $p : \R \to [0, \infty)$ with bounded support, $\mathrm{spt}(p) = [L, R]$, where $L > - \infty$ and $R < \infty$. The proofs of our previous characterizations rely on the fact that $\lim_{x \, \to \, \pm \infty} p(x) = 0$. However, we can do more: The results can be extended to cases where the limit to one endpoint of the support merely exists. The techniques needed for the proofs of the statements in this section resemble the ones we have used so far, so we shorten the arguments. As in Section \ref{SEC semi-bounded support}, we start with the characterizations via the density function before deriving further results from them. We divide the study into density functions for which the limit to the right endpoint of the support exists and such density functions for which the limit to the left endpoint exists.
\begin{lemma} \label{LEMMA bounded support, chara. via density function, right limit exists}
	Let $p$ be a probability density function with $\mathrm{spt}(p) = [L, R]$ that satisfies conditions (C1) -- (C5), and for which the limit $\lim_{x \, \nearrow \, R} p(x)$ exists. Take $X : \Omega \to (L, R)$ to be a random variable with density function $\mathfrak{f}_X$, and $\E \left| \tfrac{p^{\prime}(X)}{p(X)} \right| < \infty$. Then $X \sim p \mathcal{L}^1$ if, and only if,
	\begin{align*}
		\mathfrak{f}_X(t)
		= \mathbb{E} \left[ - \frac{p^{\prime}(X)}{p(X)} \, \mathds{1}\{ X > t \} \right] + \lim_{x \, \nearrow \, R} p(x), \quad L < t < R.
	\end{align*}
\end{lemma}

The main ideas of the proof are summarized in Appendix \ref{APP proof bounded supp chara}.

\begin{remark} \label{REMARK condition (C3) for bounded support case}
	Note that condition (C3) is simply $\int_{S(p)} | p^\prime(x) | \, \mathrm{d}x < \infty$ by the boundedness of the support. Also notice that
	\begin{align*}
		\E \left| \frac{p^\prime(X)}{p(X)} \, X \right|
		\leq \max\big\{ |L|, |R| \big\} \, \E \left| \frac{p^\prime(X)}{p(X)} \right|,
	\end{align*}
	so we never have to state both integrability conditions on $X$.
\end{remark}

\begin{remark} \label{REMARK alternative condition (C3) bounded support case, density chara.}
	By the argument of Remark \ref{REMARK integr. condition in sufficiency part of density char. semi-bounded sup.}, in the case of a continuously differentiable density with $L = 0$, we can replace the integrability condition on $X$ completely with $\E \left| \tfrac{p^\prime(X)}{p(X)} \, X \right| < \infty$, and substitute (C3) with $\int_{0}^R x | p^\prime(x) | \, \mathrm{d}x < \infty$, which is weaker than $\int_0^R |p^\prime(x)| \, \mathrm{d}x < \infty$. However, the equality
	\begin{align*}
		\mathfrak{f}_X(t)
		= \mathbb{E} \left[ - \frac{p^{\prime}(X)}{p(X)} \, \mathds{1}\{ X > t \} \right] + \lim_{x \, \nearrow \, R} p(x)
	\end{align*}
	in Lemma \ref{LEMMA bounded support, chara. via density function, right limit exists} will then only hold for $\mathcal{L}^1$-a.e. $0 < t < R$.
\end{remark}

Complementary to Lemma \ref{LEMMA bounded support, chara. via density function, right limit exists} (and with a similar proof), we have the following result.
\begin{lemma} \label{LEMMA bounded support, chara. via density function, left limit exists}
	Let $p$ be a probability density function with $\mathrm{spt}(p) = [L, R]$ that satisfies the conditions (C1) -- (C5), and for which the limit $\lim_{x \, \searrow \, L} p(x)$ exists. Let $X : \Omega \to (L, R)$ be a random variable with density function $\mathfrak{f}_X$, and $\E \left| \tfrac{p^{\prime}(X)}{p(X)} \right| < \infty$. Then $X \sim p \mathcal{L}^1$ if, and only if,
	\begin{align*}
		\mathfrak{f}_X(t)
		= \mathbb{E} \left[ \frac{p^{\prime}(X)}{p(X)} \, \mathds{1}\{ X \leq t \} \right] + \lim_{x \, \searrow \, L} p(x), \quad L < t < R.
	\end{align*}
\end{lemma}
\vspace{3mm}

With obvious adaptations, Remark \ref{REMARK alternative condition (C3) bounded support case, density chara.} also applies here (in the case $R = 0$). We now use the Lemmata \ref{LEMMA bounded support, chara. via density function, right limit exists} and \ref{LEMMA bounded support, chara. via density function, left limit exists} to derive the corresponding characterization results via the distribution function. We start again with the case of an existing limit to the right endpoint of the support.
\begin{corollary} \label{COR bounded support via distribution function, limit to right endpoint exists}
	Let $p$ be a probability density function with $\mathrm{spt}(p) = [L, R]$ such that (C1) -- (C5) are satisfied. Assume that the limit $\lim_{x \, \nearrow \, R} p(x)$ exists. Suppose that $X : \Omega \to (L, R)$ is a random variable with $\mathbb{P}\big( X \in S(p) \big) = 1$, and $\E \left| \frac{p^{\prime}(X)}{p(X)} \right| < \infty$. Then $X \sim p \mathcal{L}^1$ if, and only if,
	\begin{align*}
		F_X(t) = \mathbb{E} \left[- \frac{p^{\prime}(X)}{p(X)} \Big( \min\{ X, t \} - L \Big) \right] + (t - L) \lim_{x \, \nearrow \, R} p(x) , \quad L < t < R.
	\end{align*}
\end{corollary}
The proof runs along the lines of Theorem \ref{THM chara distr. function, sup. bounded from below}.

\begin{remark} \label{REMARK alternative condition (C3) bounded support case, distribution function chara.}
	Whenever $p$ is continuously differentiable with support $[0, R]$, it suffices to have $\int_0^R x |p^{\prime}(x)| \, \mathrm{d}x < \infty$, instead of (C3), and the weaker condition $\E\left| \frac{p^{\prime}(X)}{p(X)} \, X \right| < \infty$ to cover the requirements of Corollary \ref{COR bounded support via distribution function, limit to right endpoint exists}.
\end{remark}

The following result is complementary to Corollary \ref{COR bounded support via distribution function, limit to right endpoint exists}.
\begin{corollary} \label{COR bounded support via distribution function, limit to left endpoint exists}
	Assume that $p$ is a probability density function with $\mathrm{spt}(p) = [L, R]$ that satisfies (C1) -- (C5). Further suppose that the limit $\lim_{x \, \searrow \, L} p(x)$ exists. Let $X : \Omega \to (L, R)$ be a random variable with $\mathbb{P}\big( X \in S(p) \big) = 1$, and $\E \left| \tfrac{p^{\prime}(X)}{p(X)} \right| < \infty$. Then $X \sim p \mathcal{L}^1$ if, and only if, the distribution function of $X$ satisfies
	\begin{align*}
		1 - F_X(t) = \mathbb{E} \left[\frac{p^{\prime}(X)}{p(X)} \Big( R - \max\{ X, t \} \Big) \right] + (R - t) \lim_{x \, \searrow \, L} p(x) , \quad L < t < R.
	\end{align*}
\end{corollary}

Remark \ref{REMARK alternative condition (C3) bounded support case, distribution function chara.} applies, with minor (but obvious) adaptations, in the case $R = 0$. In general, the characterization results for probability density functions with bounded support give a good handle on a variety of wrapped and truncated distributions, like the truncated normal- or the wrapped exponential distribution. However, we state only the uniform- and the beta distribution as examples explicitly. Again, we refrain from giving the details of the calculations to check the regularity conditions. For the beta distribution, we invoke Remark \ref{REMARK alternative condition (C3) bounded support case, distribution function chara.}.
\begin{example}[\textbf{Uniform distribution}] \label{EXAMPLE Uniform}
	For $x \in (L, R)$ let $p(x) = \tfrac{1}{R - L}$ be the density function of the uniform distribution on the interval $(L, R)$. The conditions (C1) -- (C5) are trivial to check. Since the derivate of $p$ vanishes on $(L, R)$, the identities from the Corollaries \ref{COR bounded support via distribution function, limit to right endpoint exists} and \ref{COR bounded support via distribution function, limit to left endpoint exists} are the same. They read as follows. A random variable $X : \Omega \to (L, R)$ is distributed uniformly over $(L, R)$ if, and only if, its distribution function has the form
	\begin{align*}
		F_X(t) = \frac{t - L}{R - L}, \quad L < x < R.
	\end{align*}
	Apparently, we recovered the observation that the explicitly calculable form of the uniform distribution function uniquely identifies this distribution, so our characterization is redundant in this case.
\end{example}

\begin{example}[\textbf{Beta distribution}] \label{EXAMPLE Beta}
	Let $\alpha > 0$, $\beta > 1$, and
	\begin{align*}
		p(x) = \frac{x^{\alpha - 1} \, (1 - x)^{\beta - 1}}{B(\alpha, \beta)}, \quad 0 < x < 1,
	\end{align*}
	where $B(\alpha, \beta) = \frac{\Gamma(\alpha) \, \Gamma(\beta)}{\Gamma(\alpha + \beta)}$ denotes the Beta function. Since $\beta > 1$, the limit to the right endpoint of the support exists. More precisely, we have $\lim_{x \, \nearrow \, 1} p(x) = 0$. Therefore, Corollary \ref{COR bounded support via distribution function, limit to right endpoint exists} yields the following characterization. Suppose $X$ is a random variable which takes values in $(0, 1)$ almost surely and satisfies $\E \left| \tfrac{X}{1 - X} \right| < \infty$. Then $X$ has the Beta distribution with parameters $\alpha > 0$ and $\beta > 1$ if, and only if, the distribution function of $X$ has the form
	\begin{align*}
		F_X(t)
		= \E \left[ \left( \frac{\beta - 1}{1 - X} - \frac{\alpha - 1}{X} \right) \min\{ X, t \} \right], \quad 0 < t < 1.
	\end{align*}
	The Beta distribution also marks a limitation of our characterizations. If $0 < \alpha, \beta < 1$, our results fail to hold since none of the required limits exist. A special case for this phenomenon is the Arcsine distribution, which is the Beta distribution with parameters $\alpha, \beta = \tfrac{1}{2}$.
\end{example}

\section{Applications to goodness-of-fit testing}
\label{goodness-of-fit}
The idea to use distributional characterizations as a basis for statistics in goodness-of-fit testing problems is classic, see \cite{N:2017} and \cite{O'RS:1982}. In this spirit and regarding the results of the previous sections, we propose goodness-of-fit tests for any distribution with a density function that satisfies the regularity conditions of either of our characterizations (Theorems \ref{THM chara on the real line}, \ref{THM chara distr. function, sup. bounded from below}, and Corollaries \ref{COR chara support bounded from above via distribution function}, \ref{COR bounded support via distribution function, limit to right endpoint exists}, \ref{COR bounded support via distribution function, limit to left endpoint exists}). For the sake of readability, we give the following discussion in the case of continuously differentiable and positive density functions on the positive axis, dealt with in Corollary \ref{COR chara supported by pos. half axis, via distribution function}. This case includes the largest class of examples we gave previously, and it also includes the new test we provide at the end of this section. The arguments for using the characterizations for density functions on the whole real line or such densities that have bounded support to construct corresponding goodness-of-fit tests are very similar, of course.

We consider a parametric family of distributions $\mathcal{P}_\Theta = \left\{ p_\vartheta \mathcal{L}^1 \, | \,  \vartheta \in \Theta \right\}$, $\Theta\subset\R^d$, where we assume that $\mathrm{spt}(p_{\vartheta}) = [0, \infty)$ and that $p_{\vartheta}$ is continuously differentiable and positive on $(0, \infty)$. Moreover, $p_{\vartheta}$ is taken to satisfy the prerequisites of Corollary \ref{COR chara supported by pos. half axis, via distribution function}. Testing the fit of a positive random variable $X$ to $\mathcal{P}_\Theta$ means to test the hypothesis
\begin{equation} \label{hypothesis H_0}
	\mathbf{H_0}: \, \mathbb{P}^{X} \in \mathcal{P}_\Theta
\end{equation}
against general alternatives. Let $s : (0, \infty) \times \Theta \to (0, \infty)$ be a measurable function, used for scaling, such that $X \sim p_\vartheta \mathcal{L}^1$ if, and only if, $s(X; \vartheta) \sim p_{\vartheta^{*}} \mathcal{L}^1$ for some $\vartheta^{*} \in \Theta^{*} \subset \Theta$. We assume that
\begin{equation*}
	\E \left| \frac{p^{\prime}_{\vartheta^*}\big(s(X; \vartheta)\big)}{p_{\vartheta^*}\big(s(X; \vartheta)\big)} \, s(X; \vartheta) \right| < \infty.
\end{equation*}
By Corollary \ref{COR chara supported by pos. half axis, via distribution function}, we have $s(X; \vartheta) \sim p_{\vartheta^*} \mathcal{L}^1$ if, and only if, the distribution function of $s(X; \vartheta)$ has the form
\begin{align} \label{chara. for GoF}
	F_{s(X; \vartheta)}(t)
	= \E \left[ - \frac{p^{\prime}_{\vartheta^*}\big(s(X; \vartheta)\big)}{p_{\vartheta^*}\big(s(X; \vartheta)\big)} \, \min\Big\{ s(X; \vartheta), \, t \Big\} \right], \quad t > 0.
\end{align}
In order to test $\mathbf{H_0}$ based on a sample $X_1, \dots, X_n$ of independent and identically distributed (iid.) positive random variables, put $Y_{n, j} = s(X_j; \widehat{\vartheta}_n)$, $j = 1, \dots, n$. We consider the empirical distribution function $\widehat{F}_n$ of $Y_{n, 1}, \dots, Y_{n, n}$ as an estimator of $F_{s(X; \vartheta)}$. Hereby denoting a consistent estimator of $\vartheta$ by $\widehat{\vartheta}_n = \widehat{\vartheta}_n(X_1, \dots, X_n)$, we use $\widehat{\vartheta}^{*}_n = \widehat{\vartheta}_n (Y_{n, 1}, \dots, Y_{n, n})$ as an estimator of $\vartheta^{*} \in \Theta^{*}$, and take
\begin{align*}
	\widehat{T}_n(t)
	= - \frac{1}{n} \sum\limits_{j = 1}^{n} \frac{p^{\prime}_{\widehat{\vartheta}^{*}_n}\big(Y_{n , j}\big)}{p_{\widehat{\vartheta}^{*}_n}\big(Y_{n , j}\big)} \, \min\big\{ Y_{n , j}, t \big\} , \quad t > 0,
\end{align*}
as an estimator of the second quantity in (\ref{chara. for GoF}). Taking some metric $\delta$ on a set containing both functions, we propose as a goodness-of-fit statistic the quantity
\begin{align*}
	\delta \Big( \widehat{T}_n , \widehat{F}_n \Big).
\end{align*}
By (\ref{chara. for GoF}), this term ought to be close to zero under $\mathbf{H_0}$, so large values of the statistic will lead us to reject the hypothesis.

As witnessed by \cite{BH:2000}, \cite{BE:2019:2}, and \cite{BE:2019:1}, tests of this type are noteworthy competitors to established tests. An advantage lies in the range of their applicability. A substantial proportion of known procedures relies on a comparison between theoretical moment generating functions, see \cite{CQ:2005}, \cite{HJG:2018}, and \cite{Z:2010}, or characteristic functions, as employed by \cite{BH:1988}, \cite{EP:1983}, and \cite{ACJM:2009}, and their empirical pendants, or on a differential equation that characterizes the Laplace transformation, see \cite{HK:2002} and \cite{HME:2012}. All of these share the unpleasant feature that in order to establish the theoretic basis for the test statistics, one has to have explicit knowledge about these transformations for the distribution in consideration. Since their handling is not possible for every distribution, our suggestions provide a genuine alternative, for they require no more than the knowledge of the density function and its derivative. Moreover, our tests do not rely on a characterization that is tailored to one specific distribution. Instead, we provide a framework for testing fit to many different distributions, as indicated by our list of examples.

\subsection{Tests for the Gamma distribution}
\label{subsection Gamma tests}
In \cite{BE:2019:2}, the authors establish the result of Corollary \ref{COR chara supported by pos. half axis, via distribution function} for the special case of the Gamma distribution and examine the corresponding goodness-of-fit statistic. Denote by $p_{\vartheta}(x) = \tfrac{\lambda^{-k}}{\Gamma(k)} \, x^{k - 1} \, e^{-x/\lambda}$, $x > 0$, where $\vartheta = (k, \lambda) \in (0, \infty)^2 = \Theta$, the density function of the Gamma distribution with shape parameter $k$ and scale parameter $\lambda$. Let $X$ be a positive random variable with $\E X < \infty$. To reflect the scale invariance of the class of Gamma distributions, choose the scaling function $s(x; \vartheta) = x / \lambda$. Apparently, $X \sim p_{\vartheta} \mathcal{L}^1$ if, and only if, $s(X; \vartheta) \sim p_{\vartheta^*} \mathcal{L}^1$, where $\vartheta^* = (k, 1) \in (0, \infty) \times \{1\} = \Theta^*$, and
\begin{align*}
	\E \left| \frac{p^{\prime}_{\vartheta^*}\big(s(X; \vartheta)\big)}{p_{\vartheta^*}\big(s(X; \vartheta)\big)} \, s(X; \vartheta) \right|
	\leq |k - 1| + \lambda^{-1} \E X
	< \infty.
\end{align*}
By Example \ref{EXAMPLE Gamma}, $X$ follows a Gamma law with parameter vector $\vartheta = (k, \lambda)$ if, and only if,
\begin{align*}
	F_{X / \lambda}(t)
	= F_{s(X; \vartheta)}(t)
	= \E \left[ \left( - \frac{k - 1}{s(X; \vartheta)} + 1 \right) \min\big\{ s(X; \vartheta), t \big\} \right], \quad t > 0.
\end{align*}
To construct the goodness-of-fit test, let $X_1, \dots, X_n$ be iid. copies of $X$ and consider a consistent, scale equivariant estimator $\widehat{\lambda}_n = \widehat{\lambda}_n(X_1 \dots, X_n)$ of $\lambda$ as well as a consistent, scale invariant estimator $\widehat{k}_n = \widehat{k}_n(X_1, \dots, X_n)$ of $k$. We set $Y_{n, j} = s(X_j; \widehat{k}_n, \widehat{\lambda}_n) = X_j / \widehat{\lambda}_n$,  $j = 1, \dotso, n$. Naturally, $\widehat{\lambda}_n^* = \widehat{\lambda}_n(Y_{n, 1}, \dots, Y_{n, n}) = 1$ and $\widehat{k}_n^* = \widehat{k}_n(Y_{n, 1}, \dots, Y_{n, n}) = \widehat{k}_n(X_1, \dots, X_n) = \widehat{k}_n$ are consistent estimators of $\lambda^* = 1$ and $k^* = k$. In accordance with our general consideration above, we take
\begin{align*}
	\widehat{T}_n(t) = \frac{1}{n} \sum\limits_{j=1}^{n} \left( - \frac{\widehat{k}_n - 1}{Y_{n ,j}} + 1 \right) \min\{ Y_{n, j}, t \}, \quad t > 0.
\end{align*}
\cite{BE:2019:2} considered the functions $\widehat{T}_n$ and $\widehat{F}_n = n^{-1} \sum_{j=1}^{n} \mathds{1}\{ Y_{n, j} \leq \cdot \}$ as random elements of the Hilbert space $L^2\big( (0, \infty), \, \mathcal{B}^1_{>0}, \, w(t) \, \mathrm{d}t \big)$, where $w$ is an appropriate weight function. They obtained the statistic
\begin{align*}
	G_n = n \int_0^{\infty} \left| \widehat{T}_n(t) - \widehat{F}_n(t) \right|^2 w(t) \, \mathrm{d}t,
\end{align*}
derived the limit distribution under the hypothesis using the Hilbert space central limit theorem, and gave a proof of the consistency of this test procedure against fixed alternatives with existing expectation. Moreover, they explained how to implement the test using a parametric bootstrap and showed in a Monte Carlo simulation study that the test excels classical procedures and keeps up with the best Gamma tests proposed so far. Contributions like \cite{HME:2012}, \cite{PS:2017}, and \cite{VGE:2015} indicate that testing fit to the Gamma distribution is also a topic of ongoing research.

The characterization of the exponential distribution via the mean residual life function is a special case of Corollary \ref{COR chara supported by pos. half axis, via distribution function} (cf. Example \ref{EXAMPLE Exponential}), and thus the corresponding test for exponentiality is to be seen as a special case of the test for the Gamma distribution at hand. \cite{BH:2000} used the characterization, which was known in a different disguise already, to construct the associated test for exponentiality in the sense described above. They showed that the limit distribution under the hypothesis coincides with the limiting null distribution of the classical Cram\'{e}r-von Mises statistic when testing for uniformity over the unit interval. Furthermore, they proved the consistency of the test procedure against any fixed alternative distribution. The test has already been included in the extensive comparative simulation study conducted by \cite{ASSV:2017}. Adding a tuning parameter to the weight function leads to the test statistic proposed by \cite{BH:2008}. The recent papers by \cite{CMO:2018}, \cite{JMNO:2015}, \cite{N:2017}, \cite{N:2015}, \cite{TMG:2018}, \cite{VN:2015}, and \cite{ZPM:2015} show that tests for exponentiality are still of importance to the research community.

\subsection{Tests for normality}
\label{subsection normality tests}
The goodness-of-fit tests for normality proposed by \cite{BE:2019:1} are also included in our framework (cf. Example \ref{EXAMPLE Gaussian}). To fix notation, we write $p_{\vartheta}$ for the normal distribution density with mean-variance-parameter vector $\vartheta = (\mu, \sigma^2) \in \R \times (0, \infty) = \Theta$. Consider a real-valued random variable $X$ with $\E X^2 < \infty$. Taking into account the invariance under linear transformations of the class of normal distributions, \cite{BE:2019:1} used the scaling function $s(x; \vartheta) = (x - \mu) / \sigma$. Naturally, $X \sim p_{\vartheta} \mathcal{L}^1$ if, and only if, $s(X; \vartheta) \sim p_{\vartheta^*} \mathcal{L}^1$, where $\vartheta^* = (0, 1)$, that is, if $s(X; \vartheta)$ follows the standard Gaussian law. Furthermore, we have
\begin{align*}
	\E \left| \frac{p^{\prime}_{\vartheta^*}\big(s(X; \vartheta)\big)}{p_{\vartheta^*}\big(s(X; \vartheta)\big)} \right|
	= \E \big| s(X; \vartheta) \big|
	\leq \frac{1}{\sigma} \Big( \E|X| + |\mu| \Big)
	< \infty
\end{align*}
and
\begin{align*}
	\E \left| \frac{p^{\prime}_{\vartheta^*}\big(s(X; \vartheta)\big)}{p_{\vartheta^*}\big(s(X; \vartheta)\big)} \, s(X; \vartheta) \right|
	= \E \big( s(X; \vartheta) \big)^2
	\leq \frac{1}{\sigma^2} \Big( \E X^2 + 2 |\mu| \, \E|X| + \mu^2 \Big)
	< \infty.
\end{align*}
As a consequence, Example \ref{EXAMPLE Gaussian} states that $X$ follows a normal distribution with parameter vector $\vartheta = (\mu, \sigma^2)$ if, and only if,
\begin{align*}
	F_{s(X; \vartheta)}(t)
	= \mathbb{E} \Big[ s(X; \vartheta) \, \big( s(X; \vartheta) - t \big) \, \mathds{1}\{ s(X; \vartheta) \leq t \} \Big], \quad t \in \R.
\end{align*}
For iid. copies $X_1, \dots, X_n$ of $X$, we consider the sample mean $\overline{X}_n = n^{-1} \sum_{j=1}^{n} X_j$ and variance $S^2_n = n^{-1} \sum_{j=1}^{n} (X_j - \overline{X}_n)^2$ as consistent estimators of $\mu$ and $\sigma^2$. We put
\begin{align*}
	Y_{n, j} = s(X_j; \overline{X}_n, S_n^2) = (X_j - \overline{X}_n) / S_n, \quad j = 1, \dotso, n,
\end{align*}
and notice that $\widehat{\vartheta}_n^* = \big( {\overline{X}}_n^*, {S^2_n}^* \big) = (0, 1)$. Thus, we take
\begin{align*}
	\widehat{T}_n(t) = \frac{1}{n} \sum_{j=1}^{n} Y_{n, j} \, (Y_{n,j} - t) \, \mathds{1}\{ Y_{n, j} \leq t \}, \quad t \in \R.
\end{align*}
It remains to compare $\widehat{T}_n$ with the empirical distribution function $\widehat{F}_n$ of $Y_{n, 1}, \dots, Y_{n, n}$ by an appropriate measure of deviation. In particular, \cite{BE:2019:1} considered $\widehat{T}_n$ and $\widehat{F}_n$ as random elements in the Hilbert space $L^2\big( \R, \, \mathcal{B}^1, \, w(t) \, \mathrm{d}t \big)$, where $w$ is a suitable weight function, and chose as a metric the one induced by the Hilbert space norm. In accordance with our general considerations at the beginning of this section, their statistic has the form
\begin{align*}
	G_n = n \int_\R \left| \widehat{T}_n^X(t) - \widehat{F}_n(t) \right|^2 w(t) \, \mathrm{d}t.
\end{align*}
Besides specifying weight functions for which the statistic has an explicit formula, \cite{BE:2019:1} used the central limit theorem for random elements in separable Hilbert spaces to derive the limit distributions under the hypothesis $\mathbf{H_0}$ in (\ref{hypothesis H_0}). Furthermore, they established the consistency of the test procedures against fixed alternatives with existing second moment, and showed in a Monte Carlo simulation study that these tests are serious competitors to established procedures. The problem of testing for normality is still of interest in research, as evidenced by \cite{HJG:2018}, \cite{HJGM:2018}, \cite{HK:2017}, and numerous preprints.

\subsection{Classical Procedures}
\label{subsection classical procedures}
We consider the uniform distribution on the unit interval, $p(t) = \mathds{1}_{(0, 1)}(t)$, $t \in \R$. According to Example \ref{EXAMPLE Uniform}, our characterization results for the uniform distribution reduce to the fact the law is determined uniquely by its distribution function $F(t) = t$, $0 < t < 1$. Thus, in line with the general construction above, we obtain the statistics
\begin{align} \label{class procedures 1}
\begin{split}
	K_n = \sqrt{n} \, \sup\limits_{0 \, < \, t \, < \, 1} \Big| \widehat{F}_n(t) - F(t) \Big| \quad\quad \text{and} \quad\quad
	\omega_n^2 = n \int_0^1 \Big| \widehat{F}_n(t) - F(t) \Big|^2 \, \mathrm{d}F(t)
\end{split}
\end{align}
for testing the uniformity hypothesis. Here, $\widehat{F}_n$ is the empirical distribution function of $X_1, \dots, X_n$, which are iid. copies of a random variable $X$ with values in $(0, 1)$. Thus we have recovered in this special case the classical Kolmogorov-Smirnov- and Cram\'{e}r-von Mises statistics. Using a weight function in the integral statistic, we may also obtain the one from Anderson and Darling. For an account of the historical development of these classical procedures, a synoptic derivation of their limit distribution and an explanation on how to extend these tests to situations where the null hypothesis includes a whole (parametric) family of continuous distributions, as well as for further references, we recommend \cite{BCMCCdWGLRMS:2000}.

\subsection{A test for the Burr Type XII distribution}
\label{subsection Burr test}
In this subsection, we propose a new goodness-of-fit test for the two parameter Burr Type XII distribution $\mbox{Burr}_{XII}(k,c)$, $k,c>0$, based on the characterization given in Example \ref{EXAMPLE Burr}, fixing the scale parameter $\sigma=1$. The distribution is known under a variety of names, e.g. as the Singh-Maddala distribution or as the Pareto (IV) distribution, for details see \cite{KK:2003}, Section 6.2. We denote the density function corresponding to the $\mbox{Burr}_{XII}(k,c)$ distribution by $p_\vartheta(x) = c \, k \, x^{c - 1} (1 + x^c)^{- k - 1}$, $x > 0$, where $\vartheta = (k, c) \in (0, \infty)^2 = \Theta$. For iid. copies $X_1, \dots, X_n$ of $X$ define
\begin{equation*}
	\widehat{T}_n(t) = \frac{1}{n} \sum_{j = 1}^n \left( \widehat{c}_n \left( \widehat{k}_n + 1 \right) \frac{X_j^{\widehat{c}_n - 1}}{1 + X_j^{\widehat{c}_n}} - \frac{\widehat{c}_n - 1}{X_j} \right) \min\{ X_j, t \}, \quad t > 0,
\end{equation*}
in accordance with the general framework above, which leads to the family of $L^2$-type statistics
\begin{equation*}
	B_{n, a}	= n \int_0^\infty \left| \widehat{T}_n(t) - \widehat{F}_n(t) \right|^2 w_a(t) \, \mathrm{d}t.
\end{equation*}
Here $\widehat{k}_n$ and $\widehat{c}_n$ are consistent estimators of the parameters $k$ and $c$, $\widehat{F}_n$ is the empirical distribution function of $X_1, \dots, X_n$, and $w_a(t) = \exp(-at)$ is a weight function depending on a tuning parameter $a > 0$. Rejection of the hypothesis $\mathbf{H_0}$ in (\ref{hypothesis H_0}), i.e. that the data comes from the Burr Type XII family, is for large values of $B_{n, a}$. Writing $X_{(1)} \leq \dotso \leq X_{(n)}$ for
the order statistics of $X_1, \dots ,X_n$, we have after some tedious calculations
\begin{eqnarray*}
B_{n, a}&=& \frac{2}{n} \sum_{1 \leq j < \ell \leq n} \left\{ A^{[1]}_{(\ell), n} \left[ \frac{2 A^{[1]}_{(j), n}}{a^3} \, \big( 1 - e^{- a X_{(j)}} \big) + \frac{A^{[2]}_{(j), n}}{a^2} \, \big( e^{- a X_{(j)}} + e^{- a X_{(\ell)}} \big) \right.\right. \\
&&\left. \left. \phantom{\frac{2 A^{[1]}_{(j), n}}{a^3}} ~~~~~~~~~~~~~~~+ \frac{\widehat{c}_n - 2}{a^2} \, e^{- a X_{(j)}} - \frac{X_{(j)}}{a} \, e^{- a X_{(j)}} \right] + \frac{A^{[2]}_{(j), n}}{a} \, e^{- a X_{(\ell)}} \right\} \\
&&+ \frac{1}{n} \sum_{j = 1}^{n} \left\{ \big(A^{[1]}_{(j), n}\big)^2 \left( - \frac{2 X_{(j)}}{a^2} \, e^{- a X_{(j)}} - \frac{2}{a^3} \, e^{- a X_{(j)}} + \frac{2}{a^3} \right) \phantom{\frac{2 A^{[2]}_{(j), n}}{a}} \right. \\
&&\left. ~~~~~~~~~~~+ \frac{2(j - 1) \, \widehat{c}_n}{a^2} \, A^{[1]}_{(j), n} \, e^{- a X_{(j)}} + \frac{2 A^{[2]}_{(j), n}}{a} \, e^{- a X_{(j)}} \right\} \\
&&+ \frac{2 \widehat{c}_n}{a \, n} \sum_{j = 1}^{n} j \, e^{- a X_{(j)}} - \frac{1}{a \, n} \sum_{j = 1}^{n} e^{- a X_{(j)}},
\end{eqnarray*}
where
\begin{equation*}
	A^{[1]}_{(j), n} = \widehat{c}_n \, \left(\widehat{k}_n + 1\right) \, \frac{X_{(j)}^{\widehat{c}_n - 1}}{1 + X_{(j)}^{\widehat{c}_n}} - \frac{\widehat{c}_n - 1}{X_{(j)}}
	\quad \text{and} \quad
	A^{[2]}_{(j), n} = - \widehat{c}_n \, \left(\widehat{k}_n + 1\right) \, \frac{X_{(j)}^{\widehat{c}_n}}{1 + X_{(j)}^{\widehat{c}_n}},
\end{equation*}
which is an easily computable formula that avoids any numerical integration routines. In the following simulation study we show the effectiveness of this new test statistics in comparison to the classical procedures adapted for the composite hypothesis $\mathbf{H_0}$, namely the Kolmogorov-Smirnov test $K_n$, the Cram\'{e}r-von Mises test $CM$, the Anderson-Darling test $AD$, and the Watson test $WA$. Let $F(x; k,c) = 1 - (1 + x^c)^{-k}$, $x > 0$, denote the distribution function of $\mbox{Burr}_{XII}(k,c)$. The $K_n$-statistic is
\[K_n\ = \ \max\{D^+,D^-\},\]
where
\begin{eqnarray*}
D^+&=&\max_{j \, = \, 1, \dots ,n}\left(j/n-F\left(X_{(j)};\widehat{k}_n,\widehat{c}_n\right)\right), \\\ D^-&=&\max_{j \, = \, 1, \dots ,n}\left(F\left(X_{(j)};\widehat{k}_n,\widehat{c}_n\right)-(j-1)/n\right).
\end{eqnarray*}
The statistics of Cram\'{e}r-von Mises and Anderson-Darling are given by
\begin{align*}
	CM\ =\ \frac{1}{12n}+\sum_{j=1}^n \left(F\left(X_{(j)};\widehat{k}_n,\widehat{c}_n\right)-\frac{2j-1}{2n} \right)^2
\end{align*}
and
\begin{align*}
	AD = -n - \frac{1}{n} \sum_{j = 1}^n \left[ (2j - 1) \log F\Big(X_{(j)}; \widehat{k}_n, \widehat{c}_n \Big) + \big(2 (n - j) + 1\big) \log\left( 1 - F\Big(X_{(j)}; \widehat{k}_n, \widehat{c}_n \Big) \right) \right],
\end{align*}
respectively, whereas the $WA$-statistic takes the form
\begin{align*}
	WA = CM - n \left( \frac{1}{n} \sum_{j = 1}^n F\left(X_{(j)}; \widehat{k}_n, \widehat{c}_n\right) - \frac{1}{2} \right)^2.
\end{align*}
For all procedures the parameters are estimated via the maximum likelihood method, maximizing numerically the log-likelihood function, see \cite{JW:2009} and \cite{W:1983}. There are other estimation procedures available, like the maximum product of spacings method, see \cite{SG:1993}. Critical points are obtained for the classical tests, as well as for the new test, by the same parametric bootstrap procedure, as follows: For a given sample $X_1, \dots, X_n$ of size $n$, compute the estimators $\widehat{k}_n, \widehat{c}_n$ of $k$ and $c$. Conditionally on $\widehat{k}_n, \widehat{c}_n$, generate $100$ bootstrap samples of size $n$ from $\mbox{Burr}_{XII}(\widehat{k}_n, \widehat{c}_n)$. Calculate the value of the test statistic, say $B_j^*$ ($j = 1, 2, \dots, 100$), for each bootstrap sample. Obtain the critical value $p_n$ as $B_{(90)}^*$, where $B_{(j)}^*$ denote the ordered $B_j^*$--values, and reject $\mathbf{H_0}$ if $B_{n, a} = B_{n, a}(X_1, \dots, X_n) > p_n$.

The following (alternative) distributions are considered (all densities defined for $x > 0$ in dependence of a shape parameter $\theta > 0$):
\begin{enumerate}
 \item The Burr Type XII distribution $\mbox{Burr}_{XII}(k,c)$,
 \item the exponential distribution $\mbox{Exp}(\theta)$,
 \item the linear increasing failure rate law $LF(\theta)$ with density $(1 + \theta x) \exp( - x - \theta x^2 / 2)$,
 \item the half-normal distribution with density $(2 / \pi)^{1/2} \exp(- x^2 / 2)$, denoted by $HN$,
 \item the half-Cauchy distribution with density $2 / \big(\pi (1 + x^2)\big)$, denoted by $HC$,
 \item the Gompertz law $GO(\theta)$ having distribution function $1 - \exp[ \theta^{-1} (1 - e^x)]$,
 \item the inverse Gaussian law $IG(\theta)$ with density $\big(\theta / (2\pi)\big)^{1/2} x^{-3 / 2} \exp[- \theta (x - 1)^2 / (2x)]$,
 \item the Weibull distribution with density $\theta x^{\theta - 1} \exp(- x^\theta)$, denoted by $W(\theta)$,
 \item the inverse Weibull distribution with density $\theta (1 / x)^{\theta + 1} \exp[ - (1 / x)^\theta]$, denoted by $IW(\theta)$.

\end{enumerate}
All computations are performed using the statistical computing environment \texttt{R}, see \cite{R:2019}. In each scenario, we consider the sample sizes $n = 100$ and $n = 200$, and the nominal level of significance $\alpha$ is set to $0.1$. Each entry in Tables \ref{Tab.Burr} and \ref{Tab.Burr200} presents empirical rejection rates computed with $10~000$ Monte Carlo runs. The number of bootstrap samples in each run is fixed to $100$, and for the tuning parameter $a$ we consider the values $\{0.25, 0.5, 1, 3, 5, 10\}$. The best performing test for each distribution and sample size is highlighted for easy reference.
\begin{table}[h]
\centering
\begin{tabular}{l|rrrrrr|rrrr }
Alt./Test & $B_{0.25}$ & $B_{0.5}$ & $B_{1}$ & $B_{3}$ & $B_{5}$ & $B_{10}$ & $K_n$ & $CM$ & $AD$ & $WA$ \\\hline
$\mbox{Burr}_{XII}(1,1)$    & 10 & 10 & 10 & 10 & 10 & 10  & 10 & 10 & 10 & 11  \\
$\mbox{Burr}_{XII}(2,1)$    & 9  & 10 & 10 & 10 & 10 & 10  & 10 & 10 & 09 & 11  \\
$\mbox{Burr}_{XII}(4,1)$    & 6  & 7  & 8  & 11 & 10 & 10 & 11 & 10 & 10 & 11 \\
$\mbox{Burr}_{XII}(0.5,2)$  & 9  & 9  & 9  & 10 & 9 & 8  & 10 & 10 & 10 & 10 \\
$\mbox{Burr}_{XII}(2,0.5)$  & 10 & 10 & 10 & 10 & 10 & 11 & 10 & 10 & 10 & 10  \\\hline
$\mbox{Exp}(1)$             & 0  & 27 & \textbf{\textit{69}} & 55 & 44 & 42 & 51 & 61 & 66 & 53 \\
$LF(2)$                     & 0  & 0  & 2  & \textbf{\textit{77}} & 73 & 63 & 56 & 67 & 74 & 59 \\
$LF(4)$                     & 0  & 0  & 0  & 56 & \textbf{\textit{68}} & 60 & 48 & 57 & 64 & 46 \\
$HC$                        & 12 & 12 & 13 & 14 & 13 & 12 & 12 & 13 & \textbf{\textit{15}} & 14 \\
$HN$                        & 0  & 1  & 64 & 89 & 81 & 73 & 78 & 88 & \textbf{\textit{90}} & 78 \\
$GO(2)$                     & 0  & 4  & 90 & 99 & 98 & 93 & 97 & 99 & \textbf{\textit{100}} & 98 \\
$IG(0.5)$                   & 13 & 45 & 66 & 81 & \textbf{\textit{83}} & 82 & 52 & 64 & 72 & 61 \\
$IG(1.5)$                   & 2  & 6  & 22 & 40 & \textbf{\textit{48}} & 46 & 24 & 31 & 37 & 32 \\
$IG(3)$                     & 1  & 2  & 8  & 18 & \textbf{\textit{24}} & 24 & 16 & 21 & 23 & 23 \\
$W(0.5)$                    & \textbf{\textit{74}} & 70 & 60 & 32 & 23 & 22 & 52 & 61 & 65 & 53 \\
$W(1.5)$                    & 0  & 0  & 38 & \textbf{\textit{66}} & 54 & 47 & 52 & 60 & 65 & 52 \\
$W(3)$                      & 0  & 0  & 0  & 65 & \textbf{\textit{69}} & 56 & 52 & 61 & 65 & 52 \\
$IW(1)$                     & 46 & 50 & 56 & \textbf{\textit{66}} & 63 & 44 & 37 & 42 & 48 & 41 \\

\end{tabular}%
\caption{ Percentage of rejection for $10~000$ Monte Carlo repetitions ($n = 100$, $\alpha = 0.1$). }\label{Tab.Burr}
\end{table}

\begin{table}[h]
\centering
\begin{tabular}{l|rrrrrr|rrrr }
Alt./Test & $B_{0.25}$ & $B_{0.5}$ & $B_{1}$ & $B_{3}$ & $B_{5}$ & $B_{10}$ & $K_n$ & $CM$ & $AD$ & $WA$ \\\hline
$\mbox{Burr}_{XII}(1,1)$    & 10 & 10 & 10 & 10 & 10 & 10  & 10 & 10 & 10  & 10  \\
$\mbox{Burr}_{XII}(2,1)$    & 10 & 10 & 10 & 10 & 10 & 11  & 10 & 9  & 9   & 10  \\
$\mbox{Burr}_{XII}(4,1)$    & 7 & 8 &  10 &  10 &  10 &  10 & 10 & 10 & 10  & 11 \\
$\mbox{Burr}_{XII}(0.5,2)$  &  10 &  10 &  10 &  10 &  10 &  9  & 11 & 10 & 10  & 9 \\
$\mbox{Burr}_{XII}(2,0.5)$  &  10 &  10 &  9 &  10 &  9 &  10 & 10 & 10 & 10  & 10  \\\hline
$\mbox{Exp}(1)$             &  4 &  78 & \textbf{\textit{95}} &  83 &  71 &  66 & 81 & 88 & 92  & 82 \\
$LF(2)$                     &  0 &  0 &  19 &  \textbf{\textit{97}} &  95 &  89 & 85 & 92 & 95  & 88 \\
$LF(4)$                    &  0 &  0 &  0 &  89 & \textbf{\textit{91}} &  85 & 75 & 83 & 89  & 76 \\
$HC$                        &  14 &  13 &  15 &  \textbf{\textit{19}} &  18 &  15 & 15 & 17 & \textbf{\textit{19}}  & 18 \\
$HN$                        &  0 &  17 &  98 & \textbf{\textit{100}} &  98 &  94 & 97 & 99 & \textbf{\textit{100}} & 98 \\
$GO(2)$                     &  0 &  57 & \textbf{\textit{100}} & \textbf{\textit{100}} & \textbf{\textit{100}} & \textbf{\textit{100}} & \textbf{\textit{100}} & \textbf{\textit{100}} & \textbf{\textit{100}} & \textbf{\textit{100}} \\
$IG(0.5)$                   &  29 &  78 &  93 &  \textbf{\textit{99}} &  \textbf{\textit{99}} &  \textbf{\textit{99}} & 82 & 93 & 97 & 92 \\
$IG(1.5)$                   &  1 &  13 &  41 &  72 &  82 & \textbf{\textit{85}} & 44 & 53 &  66 & 56 \\
$IG(3)$                     &  0 &  2 &  13 &  34 &  48 &  \textbf{\textit{61}} & 28 & 37 & 46 & 40 \\
$W(0.5)$                    &  \textbf{\textit{97}} &  94 &  86 &  55 &  40 &  36 & 81 & 89 & 91 & 80 \\
$W(1.5)$                    &  0 &  6 &  84 &  91 &  81 &  70 & 80 & 88 & \textbf{\textit{92}} & 82 \\
$W(3)$                      &  0 &  0 &  6 &  93 &  \textbf{\textit{94}} &  88 & 80 & 89 & 92 & 82 \\
$IW(1)$                     &  74 &  78 &  83 &  \textbf{\textit{93}} &  \textbf{\textit{93}} &  85 & 63 & 72 & 78  & 69 \\

\end{tabular}%
\caption{ Percentage of rejection for $10~000$ Monte Carlo repetitions ($n = 200$, $\alpha = 0.1$). }\label{Tab.Burr200}
\end{table}
The simulation results show the (strong) dependence on the tuning parameter, but also, for an appropriate choice, the effectiveness of the new procedures, outperforming the classical procedures almost uniformly with the exception of the half-Cauchy-, half-normal- and Gompertz distribution, where the Anderson-Darling statistic is the most powerful test. Clearly, a data dependent choice for an optimal tuning parameter would be desirable. Unfortunately, there is no known procedure for this kind of test statistics, where the distribution under $\mathbf{H_0}$ depends on the true values of the parameters, but results for tests of location-scale families by \cite{AS:2015} and \cite{T:2019} give hope for new developments. A good compromise for practitioners concerning the choice of the tuning parameter is $a = 3$ in view of Tables \ref{Tab.Burr} and \ref{Tab.Burr200}.

Similar to the other procedures based on our approach, see Subsections \ref{subsection Gamma tests} and \ref{subsection normality tests}, we expect the statistics $B_{n, a}$ to converge under $\mathbf{H_0}$ to the square of the $L^2$-norm of a centered Gaussian process, and the tests to be consistent against fixed alternatives.

\section{Conclusions}
\label{section conclusions}
We devoted this work to the derivation of explicit characterizations for a large class of continuous univariate probability distributions. Our motivation was the fact that the characterization of the standard normal distribution as the unique fixed point of the zero-bias transformation reduces to an explicit formula for the distribution function of the transformed distribution. We extrapolated this formula to other distributions by applying the Stein type identity commonly used within the density approach. Research related to our characterizations concerns the study of distributional transformations, see \cite{GR:2005} and \cite{D:2017}. While these are constructed from scratch and are used to prove Stein type characterizations, we took such a Stein identity for granted and dropped the ambition to obtain distributional transformations. Thus, starting with more information and demanding less structure from the transformations, we established better accessible explicit characterization formulae. In the last section, we discussed an immediate application. We illustrated how to use the characterizations for the construction of goodness-of-fit tests. The corresponding procedures for the normal-, the exponential- and the Gamma distribution have already been investigated in the literature, and they show very promising performance. The great advantage of our approach lies in the wide range of its applicability. To confirm this last claim, we constructed the (to our best knowledge) first ever goodness-of-fit test focused on the Burr Type XII distribution.

\subsection*{Acknowledgements}
The authors would like to thank an associate editor as well as three anonymous reviewers for their comments and suggestions that led to a major improvement the paper.

\appendix
\section{Proofs}
\subsection{Proof of Lemma \ref{LEMMA density approach}}
\label{APP proof of density approach}

If $X \sim p \mathcal{L}^1$, any $f \in \mathcal{F}_p$ satisfies
\begin{align*}
	\E \left[ f^{\prime}(X) + \frac{p^{\prime}(X)}{p(X)} \, f(X) \right]
	&= \int_{S(p, f)} \big(f \cdot p\big)^{\prime}(x) \, \mathrm{d}x \\
	&= \int_L^{y^f_1} \big(f \cdot p\big)^{\prime}(x) \, \mathrm{d}x + \sum\limits_{\ell = 1}^{m} \int_{y^f_\ell}^{y^f_{\ell + 1}} \big(f \cdot p\big)^{\prime}(x) \, \mathrm{d}x + \int_{y^f_{m+1}}^R \big(f \cdot p\big)^{\prime}(x) \, \mathrm{d}x \\
	&= \lim_{x \, \nearrow \, R} f(x) \, p(x) - \lim_{x \, \searrow \, L} f(x) \, p(x) \\
	&~~~+ \sum\limits_{\ell = 1}^{m + 1} \left( \lim_{x \, \nearrow \, y^f_\ell} f(x) \, p(x) - \lim_{x \, \searrow \, y^f_\ell} f(x) \, p(x) \right) \\
	&= 0.
\end{align*}
For the converse, fix $t \in S(p) \setminus \mathrm{disc}(X)$ and define $f^p_t : (L, R) \to \R$ through
\begin{equation*}
	f^p_t(x) = \frac{1}{p(x)} \int_L^{x} \Big( \mathds{1}_{(L, t]} (s) - P(t) \Big) p(s) \, \mathrm{d}s.
\end{equation*}
The function $f^p_t$ is continuous, and
\begin{equation*}
	\lim_{x \, \nearrow \, R} f^p_t(x) \, p(x) = \int_L^{R} \Big( \mathds{1}_{(L, t]} (s) - P(t) \Big) p(s) \, \mathrm{d}s = P(t) - P(t) = 0.
\end{equation*}
Noting that $f^p_t(x) = \tfrac{1}{p(x)} \, P(x) \big( 1 - P(t) \big)$ for $x < t$, we also have $\lim_{x \, \searrow \, L} f^p_t(x) \, p(x) = 0$. With this representation of $f^p_t(x)$ for $x < t$, as well as with $f^p_t(x) = \tfrac{1}{p(x)} \big( 1 - P(x) \big) P(t)$ for $x > t$, we see that $f^p_t$ is differentiable on $S(p) \setminus \{t\} = S(p, f^p_t)$ with
\begin{align} \label{derivative of f_t^p}
	f_t^{p \, \prime}(x)
	&= - \frac{p^{\prime}(x)}{p(x)} \, f^p_t(x) + \mathds{1}_{(L, t]} (x) - P(t), \quad x \in S(p) \setminus \{t\}.
\end{align}
We get with condition (C2)
\begin{align*}
	\sup\limits_{x \, \in \, S(p) \, \setminus \, \{t\}} \left| \frac{p^{\prime}(x)}{p(x)} \, f_t^p(x) \right|
	\leq 2 \sup\limits_{x \, \in \, S(p)} \left| \frac{p^{\prime}(x) \min\{ P(x), \, 1 - P(x) \}}{p^2(x)} \right|
	= 2 \sup\limits_{x \, \in \, S(p)} \kappa_p(x)
	< \infty,
\end{align*}
and, by (\ref{derivative of f_t^p}),
\begin{align*}
	\sup\limits_{x \, \in \, S(p) \, \setminus \, \{t\}} \left| f_t^{p \, \prime}(x) \right|
	\leq \sup\limits_{x \, \in \, S(p) \, \setminus \, \{t\}} \left| \frac{p^{\prime}(x)}{p(x)} \, f_t^p(x) \right| + 2
	\leq 2 \sup\limits_{x \, \in \, S(p)} \kappa_p(x) + 2
	< \infty.
\end{align*}
Thus $f^p_t \in \mathcal{F}_p$ and $t_{f_t^p} = t \notin \mathrm{disc}(X)$. The assumption in the converse implication and (\ref{derivative of f_t^p}) yield
\begin{equation*}
	0
	= \E \left[ f_t^{p \, \prime}(X) + \frac{p^{\prime}(X)}{p(X)} \, f^p_t(X) \right]
	= \mathbb{P}(X \leq t) - P(t).
\end{equation*}
Hence $\mathbb{P}(X \leq t) = P(t)$ for all $t \in S(p) \setminus \mathrm{disc}(X)$. As $S(p) \setminus \mathrm{disc}(X)$ is dense in $(L, R)$ and $t \mapsto \mathbb{P}(X \leq t)$, $t \mapsto P(t)$ are right-continuous, the claim follows. \qed

\subsection{Proof of Theorem \ref{THM chara on the real line}}
\label{APP proof of real line chara}

Let $X \sim p \mathcal{L}^1$. By (C1) and (C3) we may use the fundamental theorem of calculus to obtain
\begin{align*}
	F_X(t) = P(t)
	&= \int_{-\infty}^{t} \left( p(y_1) + \sum\limits_{\ell = 1}^{k - 1} \Big( p(y_{\ell + 1}) - p(y_\ell) \Big) + p(s) - p(y_k) \right) \mathrm{d}s \\
	&= \int_{-\infty}^{t} \left( \int_{-\infty}^{y_1} p^{\prime}(x) \, \mathrm{d}x + \sum\limits_{\ell = 1}^{k - 1} \int_{y_\ell}^{y_{\ell + 1}} p^{\prime}(x) \, \mathrm{d}x + \int_{y_k}^{s} p^{\prime}(x) \, \mathrm{d}x \right) \mathrm{d}s,
\end{align*}
where $k = k(s)$ is the largest index in $\{ 1, \dots, m \}$ for which still $y_k < s$ [for all $s \leq y_1$ the $y_\ell$ need not be taken into account as $p$ is continuously differentiable on $(- \infty, s)$ in these cases]. Now, since $X$ has density function $p$, we easily see [still using (C1)] that
\begin{align*}
	\int_{y_\ell}^{y_{\ell + 1}} p^{\prime}(x) \, \mathrm{d}x
	= \mathbb{E} \left[ \frac{p^{\prime}(X)}{p(X)} \, \mathds{1}\{ y_\ell < X \leq y_{\ell + 1} \} \right],
\end{align*}
for $\ell \in \{ 1, \dots, k - 1 \}$, and similar representations for the other integrals give
\begin{align*}
	F_X(t) = \int_{-\infty}^{t} \mathbb{E} \left[ \frac{p^{\prime}(X)}{p(X)} \, \mathds{1}\{X \leq s\} \right] \mathrm{d}s
	= \mathbb{E} \left[ \frac{p^{\prime}(X)}{p(X)} \, (t - X) \, \mathds{1}\{ X \leq t \} \right] , \quad t \in \R,
\end{align*}
where, in the second equality, we used Fubini's theorem. That is admissible since Tonelli's theorem and (C3) imply, for each $t \in \R$,
\begin{align*}
	\int_{-\infty}^{t} \mathbb{E} \left[ \left| \frac{p^{\prime}(X)}{p(X)} \right| \, \mathds{1}\{X \leq s\} \right] \mathrm{d}s
	= \mathbb{E} \left[ \frac{| p^{\prime}(X) |}{p(X)} \, (t - X) \, \mathds{1}\{ X \leq t \} \right]
	\leq \int_{S(p)} \big|p^\prime(x)\big| \, \big(|t| + |x|\big) \, \mathrm{d}x
	< \infty.
\end{align*}
For the converse, assume that the distribution function of $X$ is given through the explicit formula in terms of $X$ as in the theorem. Putting
\begin{align*}
	d_p^X(t)
	= \E \left[ \frac{p^{\prime}(X)}{p(X)} \, \mathds{1}\{X \leq t\} \right], \quad t \in \R,
\end{align*}
condition (\ref{Integrability condition first theorem}) entails
\begin{align*}
	\E \left[ \int_{-\infty}^{t} \frac{|p^{\prime}(X)|}{p(X)} \, \mathds{1}\{X \leq s\} \, \mathrm{d}s \right]
	= \E \left[ \frac{|p^{\prime}(X)|}{p(X)} \, \big(t - X\big) \, \mathds{1}\{X \leq t\} \right]
	< \infty
\end{align*}
for every $t \in \R$. Thus, Fubini's theorem implies
\begin{align*}
	\int_{-\infty}^{t} d_p^X(s) \, \mathrm{d}s
	= \int_{-\infty}^{t} \E \left[ \frac{p^{\prime}(X)}{p(X)} \, \mathds{1}\{X \leq s\} \right] \mathrm{d}s
	= \E \left[ \frac{p^{\prime}(X)}{p(X)} \int_{-\infty}^{t} \mathds{1}\{X \leq s\} \, \mathrm{d}s \right]
	= F_X(t)
\end{align*}
for $t \in \R$. Since $F_X$ is increasing and $d_p^X$ is right-continuous, we conclude $d_p^X \geq 0$. Moreover, we infer
\begin{align*}
	\int_{\R} d_p^X (s) \, \mathrm{d}s
	= \lim\limits_{t \, \to \, \infty} \int_{-\infty}^{t} d_p^X (s) \, \mathrm{d}s
	= \lim\limits_{t \, \to \, \infty} F_X(t)
	= 1,
\end{align*}
for $F_X$ is a distribution function. Hence, $d_p^X$ is the density function of $X$. Using the first part of (\ref{Integrability condition first theorem}), dominated convergence gives
\begin{align*}
	\E \left[ \frac{p^{\prime}(X)}{p(X)} \right]
	= \lim\limits_{t \, \to \, \infty} \E \left[ \frac{p^{\prime}(X)}{p(X)} \, \mathds{1}\{X \leq t\} \right]
	= \lim\limits_{t \, \to \, \infty} d_p^X(t)
	= 0.
\end{align*}	
Therefore, we conclude that for each $f \in \mathcal{F}_p$
\begin{align*}
	\E\big[ f^{\prime}(X) \big]
	&= \int_{S(p, f)} f^{\prime}(s) \, d_p^X(s) \, \mathrm{d}s \\
	&= \int_{- \infty}^{y^f_1} f^{\prime}(s) \, \E\left[ \frac{p^{\prime}(X)}{p(X)} \, \mathds{1}\{X \leq s\} \right] \mathrm{d}s
	+ \sum\limits_{\ell = 1}^{m} \int_{y^f_\ell}^{y^f_{\ell + 1}} f^{\prime}(s) \, \E\left[ \frac{p^{\prime}(X)}{p(X)} \, \mathds{1}\{X \leq s\} \right] \mathrm{d}s \\
	&~~~+ \int_{y^f_{m + 1}}^{\infty} f^{\prime}(s) \, \E\left[ - \frac{p^{\prime}(X)}{p(X)} \, \mathds{1}\{X > s\} \right] \mathrm{d}s \\
	&= \E\left[ \frac{p^{\prime}(X)}{p(X)} \Big( f(y^f_1) - f(X) \Big) \mathds{1}\{X \leq y^f_1\} \right] \\
	&~~~+ \sum\limits_{\ell = 1}^{m} \E\left[ \frac{p^{\prime}(X)}{p(X)} \Big( f(y^f_{\ell + 1}) - f(X) \Big) \mathds{1}\{y^f_\ell < X \leq y^f_{\ell + 1}\} \right] \\
	&~~~+ \sum\limits_{\ell = 1}^{m} \E\left[ \frac{p^{\prime}(X)}{p(X)} \Big( f(y^f_{\ell + 1}) - f(y^f_\ell) \Big) \mathds{1}\{X \leq y^f_{\ell}\} \right] \\
	&~~~+ \E\left[ -\frac{p^{\prime}(X)}{p(X)} \Big( f(X) - f(y^f_{m + 1}) \Big) \mathds{1}\{X > y^f_{m + 1}\} \right] \\
	&= \E\left[ -\frac{p^{\prime}(X)}{p(X)} \, f(X) \right].
\end{align*}
In the third equality, Fubini's theorem is applicable since $f^{\prime}$ is bounded on $S(p, f)$ and we have (\ref{Integrability condition first theorem}). Lemma \ref{LEMMA density approach} yields the claim. \qed

\subsection{Proof of Theorem \ref{THM chara distr. function, sup. bounded from below}}
\label{APP proof semi-bounded supp chara}

Let $X \sim p \mathcal{L}^1$. By Theorem \ref{THM chara semi-bounded support via density function}, we have
\begin{align*}
	F_X(t) = P(t)
	= \int_{L}^{t} \mathbb{E} \left[ - \frac{p^{\prime}(X)}{p(X)} \, \mathds{1}\{X > s\} \right] \mathrm{d}s
	= \mathbb{E} \left[ - \frac{p^{\prime}(X)}{p(X)} \Big( \min\{ X, t \} - L \Big) \right] , \quad t > L,
\end{align*}
where Fubini's theorem is applicable since Tonelli's theorem and (C3) imply
\begin{align*}
	\int_{L}^{\infty} \mathbb{E} \left[ \left| \frac{p^{\prime}(X)}{p(X)} \right| \, \mathds{1}\{X > s\} \right] \mathrm{d}s
	&= \mathbb{E} \left[ \frac{| p^{\prime}(X) |}{p(X)} \, \big( X - L \big) \right] \\
	&\leq \int_{S(p)} |x| \, \big|p^\prime(x)\big| \, \mathrm{d}x + |L| \int_{S(p)} \big|p^\prime(x)\big| \, \mathrm{d}x \\
	&< \infty.
\end{align*}
For the converse implication, we put
\begin{align*}
	d_p^X(s)
	= \E \left[ - \frac{p^{\prime}(X)}{p(X)} \, \mathds{1}\{X > s\} \right], \quad s > L,
\end{align*}
and notice that the integrability conditions on $X$ imply
\begin{align} \label{Integrability of d_p^X in suff. part of sup. bounded from below}
	\E \left[ \int_{L}^{\infty} \frac{|p^{\prime}(X)|}{p(X)} \, \mathds{1}\{X > s\} \, \mathrm{d}s \right]
	\leq \E \left| \frac{p^{\prime}(X)}{p(X)} \, X \right| + |L| \cdot \E \left| \frac{p^{\prime}(X)}{p(X)} \right|
	< \infty.
\end{align}
Thus, Fubini's theorem gives
\begin{align*}
	\int_{L}^{t} d_p^X(s) \, \mathrm{d}s
	= \E \left[ - \frac{p^{\prime}(X)}{p(X)} \int_{L}^{t} \mathds{1}\{X > s\} \, \mathrm{d}s \right]
	= F_X(t), \quad t > L.
\end{align*}
Since $d_p^X$ is integrable by (\ref{Integrability of d_p^X in suff. part of sup. bounded from below}), dominated convergence implies that $F_X$ is continuous. Moreover, Lebesgue's differentiation theorem [see Theorem 3.21 from \cite{F:1999}, with nicely shrinking sets $E_h = (t, t + h)$, $h > 0$] implies
\begin{align*}
	d_p^X(t)
	= \lim_{h \, \searrow \, 0} \frac{1}{h} \int_{t}^{t + h} d_p^X(s) \, \mathrm{d}s
	= \lim_{h \, \searrow \, 0} \frac{F_X(t + h) - F_X(t)}{h}
	\geq 0
\end{align*}
for $\mathcal{L}^1$-a.e. $t > L$, where we used that $F_X$ is increasing. Finally,
\begin{align*}
	\int_L^\infty d_p^X(s) \, \mathrm{d}s
	= \lim_{t \, \to \, \infty} \int_L^t d_p^X(s) \, \mathrm{d}s
	= \lim_{t \, \to \, \infty} F_X(t)
	= 1.
\end{align*}
We conclude that $d_p^X$ is the density function of $X$. The claim follows from Theorem \ref{THM chara semi-bounded support via density function}. \qed

\begin{remark}
	Note that we could have proven the theorem with the same argument we used in Theorem \ref{THM chara on the real line}, since the first integrability condition on $X$ ensures that $d_p^X$ is left-continuous. However, in Remark \ref{REMARK alternative condition (C3) distr. chara semibounded case} we extended the argument of Remark \ref{REMARK integr. condition in sufficiency part of density char. semi-bounded sup.} dropping that first integrability condition in the case $L = 0$. Then we can no longer conclude the left-continuity, so we had to use the different argument via Lebesgue's differentiation theorem.
\end{remark}

\subsection{Proof of Lemma \ref{LEMMA bounded support, chara. via density function, right limit exists}}
\label{APP proof bounded supp chara}

The necessity part follows with a simple rewriting of the density function, as before. For the converse implication, assume that $X$ is as in the statement of the lemma, and that
\begin{align*}
	d_p^X(t)
	= \E \left[- \frac{p^{\prime}(X)}{p(X)} \, \mathds{1}\{X > t\} \right] + \lim_{x \, \nearrow \, R} p(x), \quad L < t < R,
\end{align*}
is the density function of $X$. Since we assume both (C4) and (C5), we have by Remark \ref{REMARK on class F_p part 2} for any $f \in \mathcal{F}_p$ [note that $f$ is continuous on $(L, R)$]
\begin{align*}
	\int_{S(p, f)} f^\prime(x) \, \mathrm{d}x
	&= \int_L^{y^f_1} f^\prime(x) \, \mathrm{d}x + \sum_{\ell = 1}^{m} \int_{y^f_\ell}^{y^f_{\ell + 1}} f^\prime(x) \, \mathrm{d}x + \int_{y^f_{m + 1}}^R f^\prime(x) \, \mathrm{d}x \\
	&= \lim_{x \, \nearrow \, y^f_1} f(x) - \lim_{x \, \searrow \, L} f(x) + \sum_{\ell = 1}^{m} \left( \lim_{x \, \nearrow \, y^f_{\ell + 1}} f(x) - \lim_{x \, \searrow \, y^f_\ell} f(x) \right) \\
	&~~~+ \lim_{x \, \nearrow \, R} f(x) - \lim_{x \, \searrow \, y^f_{m + 1}} f(x) \\
	&= 0,
\end{align*}
where the integral exists by the boundedness of $f^\prime$ and the fact that $S(p, f) \subset S(p) \subset (L, R)$ which is a bounded interval. Using this fact, the proof is concluded via Lemma \ref{LEMMA density approach} with a similar calculation as in previous proofs. \qed

\bibliography{lit-districhara}   

\begin{thebibliography}{}

\bibitem[Allison and Santana, 2015]{AS:2015}
Allison, J.~S. and Santana, L. (2015).
\newblock On a data-dependent choice of the tuning parameter appearing in
  certain goodness-of-fit tests.
\newblock {\em Journal of Statistical Computation and Simulation},
  85(16):3276--3288.

\bibitem[Allison et~al., 2017]{ASSV:2017}
Allison, J.~S., Santana, L., Smit, N., and Visagie, I. J.~H. (2017).
\newblock An `apples to apples' comparison of various tests for exponentiality.
\newblock {\em Computational Statistics}, 32(4):1241--1283.

\bibitem[Anastasiou, 2018]{A:2018}
Anastasiou, A. (2018).
\newblock Assessing the multivariate normal approximation of the maximum
  likelihood estimator from high-dimensional, heterogeneous data.
\newblock {\em Electronic Journal of Statistics}, 12(2):3794--3828.

\bibitem[Anastasiou and Gaunt, 2019]{AG:2018}
Anastasiou, A. and Gaunt, R. (2019+).
\newblock Multivariate normal approximation of the maximum likelihood estimator
  via the delta method.
\newblock {\em To appear in Brazilian Journal of Probability and Statistics}.

\bibitem[Anastasiou and Reinert, 2017]{AR:2017}
Anastasiou, A. and Reinert, G. (2017).
\newblock Bounds for the normal approximation of the maximum likelihood
  estimator.
\newblock {\em Bernoulli}, 23(1):191--218.

\bibitem[Anastasiou and Reinert, 2018]{AR:2018}
Anastasiou, A. and Reinert, G. (2018).
\newblock Bounds for the asymptotic distribution of the likelihood ratio.
\newblock {\em ArXiv e-prints}, 1806.03666.

\bibitem[Barbour, 1982]{B:1982}
Barbour, A.~D. (1982).
\newblock Poisson convergence and random graphs.
\newblock {\em Mathematical Proceedings of the Cambridge Philosophical
  Society}, 92(2):349--359.

\bibitem[Barbour, 1990]{B:1990}
Barbour, A.~D. (1990).
\newblock Stein's method for diffusion approximations.
\newblock {\em Probability Theory and Related Fields}, 84(3):297--322.

\bibitem[Barbour et~al., 1989]{B:1989}
Barbour, A.~D., Karo\'{n}ski, M., and Ruci\'{n}ski, A. (1989).
\newblock A central limit theorem for decomposable random variables with
  applications to random graphs.
\newblock {\em Journal of Combinatorial Theory, Series B}, 47(2):125--145.

\bibitem[Baringhaus and Henze, 1988]{BH:1988}
Baringhaus, L. and Henze, N. (1988).
\newblock A consistent test for multivariate normality based on the empirical
  characteristic function.
\newblock {\em Metrika}, 35(1):339--348.

\bibitem[Baringhaus and Henze, 2000]{BH:2000}
Baringhaus, L. and Henze, N. (2000).
\newblock Tests of fit for exponentiality based on a characterization via the
  mean residual life function.
\newblock {\em Statistical Papers}, 41(2):225--236.

\bibitem[Baringhaus and Henze, 2008]{BH:2008}
Baringhaus, L. and Henze, N. (2008).
\newblock A new weighted integral goodness-of-fit statistic for exponentiality.
\newblock {\em Statistics \& Probability Letters}, 78(8):1006--1016.

\bibitem[Betsch and Ebner, 2019a]{BE:2019:2}
Betsch, S. and Ebner, B. (2019a).
\newblock A new characterization of the {G}amma distribution and associated
  goodness-of-fit tests.
\newblock {\em Metrika}, https://doi.org/10.1007/s00184-019-00708-7.

\bibitem[Betsch and Ebner, 2019b]{BE:2019:1}
Betsch, S. and Ebner, B. (2019b).
\newblock Testing normality via a distributional fixed point property in the
  {S}tein characterization.
\newblock {\em TEST}, https://doi.org/10.1007/s11749-019-00630-0.

\bibitem[Braverman and Dai, 2017]{BD:2017}
Braverman, A. and Dai, J.~G. (2017).
\newblock Stein’s method for steady-state diffusion approximations of ${M} /
  \mathit{Ph} / n + {M}$ systems.
\newblock {\em The Annals of Applied Probability}, 27(1):550--581.

\bibitem[Braverman et~al., 2016]{BDF:2016}
Braverman, A., Dai, J.~G., and Feng, J. (2016).
\newblock Stein’s method for steady-state diffusion approximations: {A}n
  introduction through the {E}rlang-{A} and {E}rlang-{C} models.
\newblock {\em Stochastic Systems}, 6(2):301--366.

\bibitem[Caba\~{n}a and Quiroz, 2005]{CQ:2005}
Caba\~{n}a, A. and Quiroz, A. (2005).
\newblock Using the empirical moment generating function in testing for the
  {W}eibull and the type {I} extreme value distributions.
\newblock {\em TEST}, 14(2):417--432.

\bibitem[Carrillo et~al., 2014]{CCDO:2014}
Carrillo, C., Cidr\'{a}s, J., D\'{i}az-Dorado, E., and Obando-Monta\~{n}o,
  A.~F. (2014).
\newblock An approach to determine the {W}eibull parameters for wind energy
  analysis: {T}he case of {G}alicia ({S}pain).
\newblock {\em Energies}, 7(4):2676--2700.

\bibitem[Chatterjee and Shao, 2011]{CS:2011}
Chatterjee, S. and Shao, Q.-M. (2011).
\newblock Nonnormal approximation by {S}tein’s method of exchangeable pairs
  with application to the {C}urie--{W}eiss model.
\newblock {\em The Annals of Applied Probability}, 21(2):464--483.

\bibitem[Chen et~al., 2011]{CGS:2011}
Chen, L. H.~Y., Goldstein, L., and Shao, Q.-M. (2011).
\newblock {\em Normal approximation by {S}tein's method}.
\newblock Springer, Berlin.

\bibitem[Chwialkowski et~al., 2016]{CSG:2016}
Chwialkowski, K., Strathmann, H., and Gretton, A. (2016).
\newblock A kernel test of goodness of fit.
\newblock In {\em Proceedings of the 33rd International Conference on Machine
  Learning - Volume 48}, ICML'16, pages 2606--2615.

\bibitem[Cupari\'{c} et~al., 2018]{CMO:2018}
Cupari\'{c}, M., Milo\v{s}evi\'{c}, B., and Obradovi\'{c}, M. (2018).
\newblock New ${L}^2$-type exponentiality tests.
\newblock {\em ArXiv e-prints}, 1809.07585.

\bibitem[del Barrio et~al., 2000]{BCMCCdWGLRMS:2000}
del Barrio, E., Cuesta-Albertos, J.~A., Matr{\'a}n, C., Cs{\"o}rg{\"o}, S.,
  Cuadras, C.~M., de~Wet, T., Gin{\'e}, E., Lockhart, R., Munk, A., and Stute,
  W. (2000).
\newblock Contributions of empirical and quantile processes to the asymptotic
  theory of goodness-of-fit tests.
\newblock {\em TEST}, 9(1):1--96.

\bibitem[D{\"o}bler, 2015]{D:2015}
D{\"o}bler, C. (2015).
\newblock Stein's method of exchangeable pairs for the {B}eta distribution and
  generalizations.
\newblock {\em Electronic Journal of Probability}, 20(109):1--34.

\bibitem[D{\"o}bler, 2017]{D:2017}
D{\"o}bler, C. (2017).
\newblock Distributional transformations without orthogonality relations.
\newblock {\em Journal of Theoretical Probability}, 30(1):85--116.

\bibitem[Epps and Pulley, 1983]{EP:1983}
Epps, T.~W. and Pulley, L.~B. (1983).
\newblock A test for normality based on the empirical characteristic function.
\newblock {\em Biometrika}, 70(3):723--726.

\bibitem[Fang, 2014]{F:2014}
Fang, X. (2014).
\newblock Discretized normal approximation by {S}tein’s method.
\newblock {\em Bernoulli}, 20(3):1404--1431.

\bibitem[Folland, 1999]{F:1999}
Folland, G.~B. (1999).
\newblock {\em Real Analysis: Modern Techniques and Their Applications (Second
  Edition)}.
\newblock Pure and Applied Mathematics. John Wiley \& Sons, Inc., New York.

\bibitem[Gaunt et~al., 2017]{GPR:2017}
Gaunt, R., Pickett, A., and Reinert, G. (2017).
\newblock Chi-square approximation by {S}tein's method with application to
  {P}earson's statistic.
\newblock {\em Annals of Applied Probability}, 27(2):720--756.

\bibitem[Goldstein and Reinert, 1997]{GR:1997}
Goldstein, L. and Reinert, G. (1997).
\newblock Stein's method and the zero bias transformation with application to
  simple random sampling.
\newblock {\em The Annals of Applied Probability}, 7(4):935--952.

\bibitem[Goldstein and Reinert, 2005]{GR:2005}
Goldstein, L. and Reinert, G. (2005).
\newblock Distributional transformations, orthogonal polynomials, and {S}tein
  characterizations.
\newblock {\em Journal of Theoretical Probability}, 18(1):237--260.

\bibitem[G\"{o}tze, 1991]{G:1991}
G\"{o}tze, F. (1991).
\newblock On the rate of convergence in the multivariate {CLT}.
\newblock {\em The Annals of Probability}, 19(2):724--739.

\bibitem[Henze and Jim{\'e}nez-Gamero, 2019]{HJG:2018}
Henze, N. and Jim{\'e}nez-Gamero, M.~D. (2019).
\newblock A new class of tests for multinormality with i.i.d. and garch data
  based on the empirical moment generating function.
\newblock {\em TEST}, 28(2):499--521.

\bibitem[Henze et~al., 2019]{HJGM:2018}
Henze, N., Jim{\'e}nez-Gamero, M.~D., and Meintanis, S.~G. (2019).
\newblock Characterizations of multinormality and corresponding tests of fit,
  including for {GARCH} models.
\newblock {\em Econometric Theory}, 35(3):510--546.

\bibitem[Henze and Klar, 2002]{HK:2002}
Henze, N. and Klar, B. (2002).
\newblock Goodness-of-fit tests for the inverse {G}aussian distribution based
  on the empirical {L}aplace transform.
\newblock {\em Annals of the Institute of Statistical Mathematics},
  54(2):425--444.

\bibitem[Henze and Koch, 2017]{HK:2017}
Henze, N. and Koch, S. (2017).
\newblock On a test of normality based on the empirical moment generating
  function.
\newblock {\em Statistical Papers}, https://doi.org/10.1007/s00362-017-0923-7.

\bibitem[Henze et~al., 2012]{HME:2012}
Henze, N., Meintanis, S.~G., and Ebner, B. (2012).
\newblock Goodness-of-fit tests for the {G}amma distribution based on the
  empirical {L}aplace transform.
\newblock {\em Communications in Statistics - Theory and Methods},
  41(9):1543--1556.

\bibitem[Hudson, 1978]{H:1978}
Hudson, H.~M. (1978).
\newblock A natural identity for exponential families with applications in
  multiparameter estimation.
\newblock {\em The Annals of Statistics}, 6(3):473--484.

\bibitem[Jalali and Watkins, 2009]{JW:2009}
Jalali, A. and Watkins, A.~J. (2009).
\newblock On maximum likelihood estimation for the two parameter {B}urr {XII}
  distribution.
\newblock {\em Communications in Statistics - Theory and Methods},
  38(11):1916--1926.

\bibitem[Jim{\'e}nez-Gamero et~al., 2009]{ACJM:2009}
Jim{\'e}nez-Gamero, M.~D., Alba-Fern{\'a}ndez, V., Mu\~{n}oz{-}Garc\'{i}a, J.,
  and Chalco-Cano, Y. (2009).
\newblock Goodness-of-fit tests based on empirical characteristic functions.
\newblock {\em Computational Statistics \& Data Analysis}, 53(12):3957--3971.

\bibitem[Jovanovi\'{c} et~al., 2015]{JMNO:2015}
Jovanovi\'{c}, M., Milo\v{s}evi\'{c}, B., Nikitin, Y.~Y., Obradovi\'{c}, M.,
  and Volkova, K.~Y. (2015).
\newblock Tests of exponentiality based on {A}rnold--{V}illasenor
  characterization and their efficiencies.
\newblock {\em Computational Statistics \& Data Analysis}, 90:100--113.

\bibitem[Kim, 2000]{K:2000}
Kim, S.-T. (2000).
\newblock A use of the {S}tein-{C}hen method in time series analysis.
\newblock {\em Journal of Applied Probability}, 37(4):1129--1136.

\bibitem[Kleiber and Kotz, 2003]{KK:2003}
Kleiber, C. and Kotz, S. (2003).
\newblock {\em Statistical Size Distributions in Economics and Actuarial
  Sciences}.
\newblock Wiley Series in Probability and Statistics. John Wiley and Sons,
  Inc., Hoboken, New Jersey.

\bibitem[Ley et~al., 2017]{LRS:2017}
Ley, C., Reinert, G., and Swan, Y. (2017).
\newblock Stein’s method for comparison of univariate distributions.
\newblock {\em Probability Surveys}, 14:1--52.

\bibitem[Ley and Swan, 2011]{LS:2011}
Ley, C. and Swan, Y. (2011).
\newblock A unified approach to {S}tein characterizations.
\newblock {\em ArXiv e-prints}, 1105.4925v3.

\bibitem[Ley and Swan, 2013a]{LS:2013:2}
Ley, C. and Swan, Y. (2013a).
\newblock Local {P}insker inequalities via {S}tein's discrete density approach.
\newblock {\em IEEE Transactions on Information Theory}, 59(9):5584--5591.

\bibitem[Ley and Swan, 2013b]{LS:2013}
Ley, C. and Swan, Y. (2013b).
\newblock Stein's density approach and information inequalities.
\newblock {\em Electronic Communications in Probability}, 18.

\bibitem[Ley and Swan, 2016]{LS:2016}
Ley, C. and Swan, Y. (2016).
\newblock Parametric {S}tein operators and variance bounds.
\newblock {\em Brazilian Journal of Probability and Statistics},
  30(2):171--195.

\bibitem[Linnik, 1962]{L:1953}
Linnik, Y.~V. (1962).
\newblock Linear forms and statistical criteria {I}, {II}.
\newblock {\em Selected Translations in Mathematical Statistics and
  Probability}, 3:1--40, 41--90. Originally published 1953 in the Ukrainian
  Mathematical Journal, Vol. 5, pp. 207--243, 247--290 (in Russian).

\bibitem[Liu et~al., 2016]{LLJ:2016}
Liu, Q., Lee, J.~D., and Jordan, M. (2016).
\newblock A kernelized {S}tein discrepancy for goodness-of-fit tests.
\newblock In {\em Proceedings of the 33rd International Conference on Machine
  Learning - Volume 48}, ICML'16, pages 276--284.

\bibitem[Nikitin, 2017]{N:2017}
Nikitin, Y.~Y. (2017).
\newblock Tests based on characterizations, and their efficiencies: A survey.
\newblock {\em Acta et Commentationes Universitatis Tartuensis de Mathematica},
  21(1):3--24.

\bibitem[Noughabi, 2015]{N:2015}
Noughabi, H.~A. (2015).
\newblock Testing exponentiality based on the likelihood ratio and power
  comparison.
\newblock {\em Annals of Data Science}, 2(2):195--204.

\bibitem[O'Reilly and Stephens, 1982]{O'RS:1982}
O'Reilly, F.~J. and Stephens, M.~A. (1982).
\newblock Characterizations and goodness of fit tests.
\newblock {\em Journal of the Royal Statistical Society. Series B
  (Methodological)}, 44(3):353--360.

\bibitem[Pek{\"o}z and R{\"o}llin, 2011]{PR:2011}
Pek{\"o}z, E.~A. and R{\"o}llin, A. (2011).
\newblock New rates for exponential approximation and the theorems of
  {R}\'{e}nyi and {Y}aglom.
\newblock {\em The Annals of Probability}, 39(2):587--608.

\bibitem[Pinelis, 2017]{P:2017}
Pinelis, I. (2017).
\newblock Optimal-order uniform and nonuniform bounds on the rate of
  convergence to normality for maximum likelihood estimators.
\newblock {\em Electronic Journal of Statistics}, 11(1):1160--1179.

\bibitem[Plubin and Siripanich, 2017]{PS:2017}
Plubin, B. and Siripanich, P. (2017).
\newblock An alternative goodness-of-fit test for a {G}amma distribution based
  on the independence property.
\newblock {\em Chiang Mai Journal of Science}, 44(3):1180--1190.

\bibitem[Prakasa~Rao, 1979]{R:1979}
Prakasa~Rao, B. L.~S. (1979).
\newblock Characterizations of distributions through some identities.
\newblock {\em Journal of Applied Probability}, 16(4):903--909.

\bibitem[Proakis and Salehi, 2008]{PS:2008}
Proakis, J.~G. and Salehi, M. (2008).
\newblock {\em Digital Communications, 5th Edition}.
\newblock McGraw-Hill, New York.

\bibitem[{R Core Team}, 2019]{R:2019}
{R Core Team} (2019).
\newblock {\em R: A Language and Environment for Statistical Computing}.
\newblock R Foundation for Statistical Computing, Vienna, Austria.

\bibitem[Reinert and R\"{o}llin, 2010]{RR:2010}
Reinert, G. and R\"{o}llin, A. (2010).
\newblock Random subgraph counts and {U}-statistics: {M}ultivariate normal
  approximation via exchangeable pairs and embedding.
\newblock {\em Journal of Applied Probability}, 47(2):378--393.

\bibitem[Rogers, 2008]{R:2008}
Rogers, G.~L. (2008).
\newblock Multiple path analysis of reflectance from turbid media.
\newblock {\em Journal of the Optical Society of America A}, 25(11):2879--2883.

\bibitem[Ross, 2011]{R:2011}
Ross, N. (2011).
\newblock Fundamentals of {S}tein’s method.
\newblock {\em Probability Surveys}, 8:210--293.

\bibitem[Shah and Gokhale, 1993]{SG:1993}
Shah, A. and Gokhale, D.~V. (1993).
\newblock On maximum product of spacings (mps) estimation for {B}urr {XII}
  distributions.
\newblock {\em Communications in Statistics - Simulation and Computation},
  22(3):615--641.

\bibitem[Singh and Maddala, 1976]{SM:1976}
Singh, S.~K. and Maddala, G.~S. (1976).
\newblock A function for size distribution of incomes.
\newblock {\em Econometrica}, 44(5):963--970.

\bibitem[Singh, 1987]{S:1987}
Singh, V.~P. (1987).
\newblock On application of the {W}eibull distribution in hydrology.
\newblock {\em Water Resources Management}, 1(1):33--43.

\bibitem[Stein, 1986]{S:1986}
Stein, C. (1986).
\newblock Approximate computation of expectations.
\newblock {\em Lecture Notes - Monograph Series}, 7, Institute of Mathematical
  Statistics.

\bibitem[Stein et~al., 2004]{SDHR:2004}
Stein, C., Diaconis, P., Holmes, S., and Reinert, G. (2004).
\newblock Use of exchangeable pairs in the analysis of simulations.
\newblock In {\em Stein's Method, edited by P. Diaconis and S. Holmes},
  volume~46 of {\em Lecture Notes -- Monograph Series}, pages 1--25, Beachwood,
  Ohio, USA. Institute of Mathematical Statistics.

\bibitem[Tenreiro, 2019]{T:2019}
Tenreiro, C. (2019).
\newblock On the automatic selection of the tuning parameter appearing in
  certain families of goodness-of-fit tests.
\newblock {\em Journal of Statistical Computation and Simulation},
  89(10):1780--1797.

\bibitem[Torabi et~al., 2018]{TMG:2018}
Torabi, H., Montazeri, N.~H., and Gran\'{e}, A. (2018).
\newblock A wide review on exponentiality tests and two competitive proposals
  with application on reliability.
\newblock {\em Journal of Statistical Computation and Simulation},
  88(1):108--139.

\bibitem[Villase\~{n}or and Gonz\'{a}lez-Estrada, 2015]{VGE:2015}
Villase\~{n}or, J.~A. and Gonz\'{a}lez-Estrada, E. (2015).
\newblock A variance ratio test of fit for {G}amma distributions.
\newblock {\em Statistics \& Probability Letters}, 96(C):281--286.

\bibitem[Volkova and Nikitin, 2015]{VN:2015}
Volkova, K.~Y. and Nikitin, Y.~Y. (2015).
\newblock Exponentiality tests based on {A}hsanullah's characterization and
  their efficiency.
\newblock {\em Journal of Mathematical Sciences}, 204(1):42--54.

\bibitem[Wingo, 1983]{W:1983}
Wingo, D.~R. (1983).
\newblock Maximum likelihood methods for fitting the {B}urr type {XII}
  distribution to life test data.
\newblock {\em Biometrical Journal}, 25(1):77--84.

\bibitem[Ying, 2017]{Y:2017}
Ying, L. (2017).
\newblock Stein's method for mean-field approximations in light and heavy
  traffic regimes.
\newblock In {\em SIGMETRICS 2017 Abstracts - Proceedings of the 2017 ACM
  SIGMETRICS / International Conference on Measurement and Modeling of Computer
  Systems}. Association for Computing Machinery, Inc.

\bibitem[Zardasht et~al., 2015]{ZPM:2015}
Zardasht, V., Parsi, S., and Mousazadeh, M. (2015).
\newblock On empirical cumulative residual entropy and a goodness-of-fit test
  for exponentiality.
\newblock {\em Statistical Papers}, 56(3):677--688.

\bibitem[Zghoul, 2010]{Z:2010}
Zghoul, A.~A. (2010).
\newblock A goodness of fit test for normality based on the empirical moment
  generating function.
\newblock {\em Communications in Statistics -- Simulation and Computation},
  39(6):1292--1304.

\end{thebibliography}
\bibliographystyle{apalike}

\end{document}